\documentclass[10pt]{article}
\usepackage{extsizes,graphicx,amssymb,amsmath,amsthm,amsfonts,mathtext, color,wrapfig,ytableau,mathtools}

\def\be{\begin{eqnarray}}
\def\ee{\end{eqnarray}}

\hoffset= - 0.8in         

\begin{document}

\thispagestyle{empty}
\baselineskip13.0pt

\ytableausetup{boxsize=5pt}

\vspace{-1ex}\hfill IITP/TH-21/21

\bigskip
\bigskip

\centerline{\Large{Non-stationary difference equation for $q$-Virasoro conformal blocks
}}

\vspace{3ex}

\centerline{\large{\emph{Sh.Shakirov\footnote{University of Geneva, Switzerland; Shamil.Shakirov@unige.ch \\
Institute for Information Transmission Problems, Moscow, Russia }}}}

\vspace{3ex}

\centerline{ABSTRACT}

\bigskip

{\footnotesize
Conformal blocks of $q,t$-deformed Virasoro and ${\cal W}$-algebras are important special functions in representation theory with applications in geometry and physics. In the Nekrasov-Shatashvili limit $t \rightarrow 1$, whenever one of the representations is degenerate then conformal block satisfies a difference equation with respect to the coordinate associated with that degenerate representation. \linebreak This is a stationary Schrodinger equation for an appropriate relativistic quantum integrable system. It is expected that \linebreak generalization to generic $t \neq 1$ is a non-stationary Schrodinger equation where $t$ parametrizes shift in time. In this paper \linebreak we make the non-stationary equation explicit for the $q,t$-Virasoro block with one degenerate and four generic Verma modules, and prove it when three modules out of five are degenerate, using occasional relation to Macdonald polynomials.
}

\section{Introduction}

The Virasoro algebra \cite{Virasoro} with central charge $c$, generators $\{ L_n | n \in {\mathbb Z} \}$ and relations

\begin{align}
[L_n, L_m] = (n-m) L_{n+m} + c \ \dfrac{n^3 - n}{12} \ \delta_{n,-m}, \ \ \ \ \ n,m \in {\mathbb Z}
\label{Virasoro}
\end{align}
\smallskip\\
is a remarkable infinite-dimensional Lie algebra which arises \cite{VirasoroGF} as the unique central extension of \linebreak infinitesimal two-dimensional conformal transformations, and as such has natural and long known \linebreak connections to two-dimensional conformal field theory \cite{BPZ} and string theory \cite{Polyakov}. Since then it made an interesting pattern of appearances in algebraic geometry, from constraining intersection numbers \cite{ModuliW, ModuliK} on the moduli space of curves over ${\mathbb C}$ to acting on equivariant cohomology \cite{MO} of the moduli space of \linebreak torsion free sheaves on ${\mathbb CP}^2$. In physics these appearances are related, respectively, to loop equations for two-dimensional quantum gravity \cite{DVV} and AGT conjecture in four-dimensional supersymmetric gauge theory \cite{AGT}. \linebreak As Virasoro algebra gets increasingly recognized as a common symmetry of these and other problems, \linebreak representation theory of eq. (\ref{Virasoro}) becomes significant as an inherent part of their study.

Representations of the Virasoro algebra are well understood with an important role played by Verma modules ${\mathcal V}_{\alpha}$ for $\alpha \in {\mathbb C}$ -- infinite-dimensional vector spaces spanned by the action of negative generators on the vector of highest weight parametrized by $\alpha$. Representation ${\mathcal V}_{\alpha}$ is irreducible for generic highest weight. \linebreak It is well-known that a coproduct can be defined \cite{VirCoproduct} so that tensor products of Verma modules acquire \linebreak structure of representations of the Virasoro algebra. This allows to introduce intertwiners (invariant tensors) \linebreak between several Verma modules, also known as \textit{conformal blocks}. Originated in conformal field theory, \linebreak conformal blocks are mathematically identical to intertwiners of the Virasoro algebra.

Modern point of view \cite{EK} is to consider matrix elements and traces of intertwiners of various algebras as \linebreak distinguished special functions generalizing the famous hypergeometric \cite{SF-Hyper} Gamma \cite{SF-Gamma} theta \cite{SF-Theta} families. \linebreak This is both historically justified, since most prominent special functions in mathematics have been tied to representations of symmetries in one way or another \cite{SF-Vile, SF-Wign} and practically effective, since many old \cite{EK-2} and new \cite{GenMac-1} special functions can be constructed in this way. In physics, the value proposition of this point of view has been suggested for both quantum gravity \cite{SF-New-1} and gauge theory \cite{SF-New-2}. Among other advantages, \linebreak this is a unifying perspective: as soon as Virasoro or another algebra appears in a particular problem, \linebreak one can expect a quantitative description in terms of the special functions associated to that algebra.

Following this point of view, conformal blocks (or, more precisely, their matrix elements and traces which are also called conformal blocks by slight abuse of terminology) are the main special functions that should be \linebreak associated to Virasoro algebra. If so, one may expect conformal blocks to possess properties similar to other special functions -- for example, admit nice integral representations akin to hypergeometric functions. Indeed, such an integral representation exists and goes back to Dotsenko and Fateev \cite{DF}:

\begin{align}
{\mathcal I}(x_1, \ldots, x_m) = \oint\limits_{\gamma_1} dz_1 \ldots \oint\limits_{\gamma_N} dz_N \ \prod\limits_{i < j} (z_i - z_j)^{2 \beta} \ \prod\limits_{i} \prod\limits_{a = 1}^{m} (z_i - x_a)^{\alpha_a}
\label{ConfBlock}
\end{align}
\smallskip\\
It is an $N$-fold integral representation for the highest weight matrix element of the tree-shaped \cite{DFcomb} intertwiner \linebreak between $(m+1)$ Verma modules ${\mathcal V}_{\alpha_1}, \ldots, {\mathcal V}_{\alpha_m}; {\mathcal V}_{\alpha_{\infty}}$ of the Virasoro algebra with $c = 13 - 6 \beta - 6/\beta$, where $x_1, \ldots, x_m; \infty$ are the fusion coordinates associated with these representations due to Virasoro coproduct \cite{VirCoproduct}. \linebreak The last Verma parameter $\alpha_{\infty}$ and all Verma parameters associated to internal edges in the tree graph do not enter the integral explicitly but are determined by the number and choice of contours $\gamma_1, \ldots, \gamma_N$ \cite{DFcomb}. For non-tree graphs, i.e. traces of intertwiners in the framework of \cite{EK}, integral representation exists as well \cite{DFloop} but will not be addressed here. We focus on tree conformal blocks.

Many generalizations of Virasoro algebra can be considered, for example the ${\cal W}$-algebras \cite{OrdinaryW}. \linebreak In this paper, we are interested in another notable generalization \cite{qtVir, qtVirReview}:

\begin{align}
\sum\limits_{\ell = 0}^{\infty} \ f_{\ell} \ \big( T_{n-\ell} T_{m+\ell} - T_{m-\ell} T_{n+\ell} \big) = \dfrac{(1-q)(1-t^{-1})}{1 - q t^{-1}} (q^{-n} t^{n} - q^{n} t^{-n}) \ \delta_{n,-m}, \ \ \ \ \ n,m \in {\mathbb Z}
\label{qtVirRelations}
\end{align}
\smallskip\\
where the coefficients are given by $f_0 + f_1 z + f_2 z^2 + \ldots = \exp\big[ \sum_{n \geq 1} n^{-1} z^n (1-q^n)(1-t^{-n})/(1 + q^n t^{-n}) \big]$. \linebreak This algebraic structure is known as the quantum or $q$-Virasoro or sometimes $q,t$-deformed Virasoro algebra. \linebreak It can be shown \cite{qtVir} that in the limit $q, t = q^{\beta} \rightarrow 1$ with fixed $\beta$ these relations reduce to Virasoro relations with $c = 13 - 6 \beta - 6 / \beta$, justifying the name deformed Virasoro algebra. In contrast with ordinary Virasoro, \linebreak this is not a Lie algebra and its representation theory is less understood: while it is straightforward \linebreak to define Verma modules, no suitable coproduct structure is known. As the only alternative, one is forced to consider bosonic (free-field) Fock space representation for Verma modules \cite{qtVir} and define \linebreak vaccuum matrix elements of products of screening charge operators as the equivalent of intertwiners \cite{qtBlock}. \linebreak In the case of ordinary Virasoro algebra, two constructions (coproduct and bosonic) agree with each other, but for generic $q,t$ only the latter is available. It gives the following $q,t$-deformed Dotsenko-Fateev integral:

\begin{align}
{\mathcal I}_{q,t}(x_1, \ldots, x_m) = \oint\limits_{\gamma_1} dz_1 \ldots \oint\limits_{\gamma_N} dz_N \ \ \prod\limits_{i < j} \ z_j^{2 \beta} \ \dfrac{\varphi(z_i/z_j)\varphi(q/t \ z_i/z_j)}{\varphi(q \ z_i/z_j) \varphi(t \ z_i/z_j)} \ \prod\limits_{i} \prod\limits_{a = 1}^{m} \ z_i^{\alpha_a} \ \dfrac{\varphi(x_a/z_i)}{\varphi(q^{\alpha_a} x_a/z_i)}
\label{qtConfBlock}
\end{align}
The integral now depends on one more parameter $q$ with $t = q^{\beta}$, and $\varphi(x) = \prod_{i \geq 0} (1 - q^i x) = (x;q)_{\infty}$ is the Faddeev-Kashaev quantum dilogarithm \cite{FK} which is closely related to the $q$-Gamma function (\ref{FKdilog}). \linebreak Basic properties of the function $\varphi(x)$ imply that in the limit $q, t = q^{\beta} \rightarrow 1$ with fixed $\beta$ we have ${\mathcal I}_{q,t} \rightarrow {\mathcal I}$. We consider eq. (\ref{qtConfBlock}) as the main special function associated with the $q$-Virasoro algebra (\ref{qtVirRelations}).

Continuing the logic of comparison with hypergeometric functions, one may expect that conformal blocks ${\mathcal I}_{q,t}$ should satisfy interesting (and sufficiently simple) equations with respect to the variables $x_1, \ldots, x_m$. This expectation is too good to be true for generic values of the parameters $\alpha_1, \ldots, \alpha_m$, however it is \linebreak believed to be valid when at least one of the parameters takes a special/degenerate value, namely $+1$ or $-\beta$. \linebreak The equations for conformal blocks with one degenerate parameter are well understood in two limiting cases: \linebreak first, for $(q,t \rightarrow 1)$ which is the case of ordinary Virasoro algebra; second, for $(t \rightarrow 1, q \ \mbox{fixed})$ which is another case known as the Nekrasov-Shatashvili limit \cite{NS}. In the first case, the equations are the \linebreak celebrated Belavin-Polyakov-Zamolodchikov equations for Virasoro conformal blocks \cite{BPZ} which can be viewed \linebreak as non-stationary Schrodinger differential equations for some quantum integrable many-body systems \cite{AG-Negu}. \linebreak In the second case, the equations have a form of stationary but finite-difference Schrodinger equations \cite{statDefect}. \linebreak Unifying both perspectives, the equations for generic $q,t$ are expected to be non-stationary and difference. \linebreak Making these equations explicit for conformal blocks ${\mathcal I}_{q,t}$ is the primary goal of this paper. For the main part, \linebreak we focus on the particular case ${\mathcal I}_{q,t}(0,1,\Lambda/x,\Lambda)$, leaving generalization to $m \neq 4$ for the future.

Equations with both non-stationary and finite-difference properties have been recently considered in \cite{Shiraishi1, Shiraishi2}. The special functions they consider are neither identical to ${\mathcal I}_{q,t}(x_1, \ldots, x_m)$ nor more general; rather, they form two distinct classes with non-empty intersection. In terms of representation theory, \linebreak in this paper we consider just matrix elements of intertwiners (tree conformal blocks) while the functions they consider should correspond to traces (non-tree a.k.a. higher genus conformal blocks). In addition, most of their consideration should correspond to the higher rank ${\cal W}$-algebras rather than just Virasoro algebra. When intersection happens our results agree with theirs: cf. eq. (11) in \cite{Shiraishi1} and eq. (\ref{NonStatEqToda}) in current paper.

\section{Conformal blocks of $q,t$-deformed Virasoro algebra}

Let $q,t \in {\mathbb C}^{*}$ be generic, assuming $q,t$ are not roots of unity and $t = q^{\beta}$ with $\beta \notin {\mathbb Q}$.

\paragraph{Definition 1.}Let ${\mathcal H}$ be the associative algebra generated by $\{ h_n \big| n \in {\mathbb Z} \}$ with relations

\begin{align}
[ h_n, h_m ] = n \ \dfrac{1 - q^{|n|}}{1 - t^{|n|}} \ \delta_{n,-m}, \ \ \ n,m \in {\mathbb Z}
\end{align}
\smallskip\\
where the generators are called boson modes. Note this is a Heisenberg algebra except the conjugate of $h_0$ is not introduced (it will appear shortly when we consider representations of ${\mathcal H}$).

\paragraph{Definition 2.}Let ${\mathcal F}_{\sigma}$ for $\sigma \in {\mathbb C}$ be the vector space

\begin{align}
{\mathcal F}_{\sigma} = {\mathbb C}[ h_{-1}, h_{-2}, \ldots ] \ \big| \sigma \big>
\end{align}
\smallskip\\
with the structure of representation of ${\mathcal H}$ given by

\begin{align}
h_n \big| \sigma \big> = \frac{\sigma}{2 \beta} \ \delta_{n,0} \ \big| \sigma \big>, \ \ \ n \geq 0
\end{align}
\smallskip\\
This is known as the Fock space for ${\mathcal H}$.

\paragraph{Definition 3.}Let $e^{\alpha Q}: {\mathcal F}_{\sigma} \rightarrow {\mathcal F}_{\sigma + 2 \alpha}$ for $\alpha \in {\mathbb C}$ be the linear operators defined by

\begin{align}
e^{\alpha Q} \big| \sigma \big> = \big| \sigma + 2 \alpha \big>, \ \ \ \sigma, \alpha \in {\mathbb C} \\
\nonumber \\
[h_n, e^{\alpha Q}] = 0, \ \ \ n \neq 0, \ \ \ \alpha \in {\mathbb C}
\end{align}
\smallskip\\
Note this would imply $[h_0, Q] = \beta^{-1}$ if $Q$ was defined on its own, i.e. in this case $Q$ would be the Heisenberg conjugate of $h_0$. However for what follows we only need operators $e^{\alpha Q}$ and not $Q$ itself.

\newpage

\paragraph{Definition 4.}Let $T_n: {\mathcal F}_{\sigma} \rightarrow {\mathcal F}_{\sigma}$ for $n \in {\mathbb Z}$ be the linear operators defined by

\begin{align}
\nonumber T(z) = \sum\limits_{n \in {\mathbb Z}} \ T_n \ z^{-n} \ = \ & v \ \exp\left( - \sum\limits_{k = 1}^{\infty} z^k v^{-k} \dfrac{1 - t^k}{q^k + t^k} \dfrac{h_{-k}}{k} \right)
\exp\left( - \sum\limits_{k = 1}^{\infty} z^{-k} v^{k} (1 - t^k) \dfrac{h_{k}}{k} \right) q^{\beta a_0} \\
& \emph{} + v^{-1} \exp\left( \sum\limits_{k = 1}^{\infty} z^k v^{k} \dfrac{1 - t^k}{q^k + t^k} \dfrac{h_{-k}}{k} \right)
\exp\left( \sum\limits_{k = 1}^{\infty} z^{-k} v^{-k} (1 - t^k) \dfrac{h_{k}}{k} \right) q^{-\beta a_0}
\end{align}
\smallskip\\
where $v = q^{\frac{1}{2}} t^{\frac{-1}{2}}$. Each operator $T_n$ is well-defined because in $T(z)$ the exponentials of positive generators $h_k, k > 0$ appear to the right of the exponentials of the negative generators. For any given vector in ${\mathcal F}_{\sigma}$ there are only finitely many operators of the form $h_{k_1} \ldots h_{k_m}$ with $k_1, \ldots, k_m > 0$ that do not annihilate this vector. Therefore the action of $T_n$ on any vector in ${\mathcal F}_{\sigma}$ is a finite sum of vectors in ${\mathcal F}_{\sigma}$, so $T_n$ is well-defined.

\paragraph{Proposition 5. \cite{qtVir}}The operators $\{ T_n \big| n \in {\mathbb Z} \}$ satisfy identities

\begin{align}
\sum\limits_{\ell = 0}^{\infty} \ f_{\ell} \ \big( T_{n-\ell} T_{m+\ell} - T_{m-\ell} T_{n+\ell} \big) = \dfrac{(1-q)(1-t^{-1})}{1 - q t^{-1}} (q^{-n} t^{n} - q^{n} t^{-n}) \ \delta_{n,-m}, \ \ \ \ \ n,m \in {\mathbb Z}
\label{qtVirRelations2}
\end{align}
\smallskip\\
where coefficients $f_{\ell}, \ell \geq 0$ are defined by

\begin{align}
f(z) = \sum\limits_{\ell = 0}^{\infty} f_{\ell} z^{\ell} = \exp\left( \sum\limits_{n = 1}^{\infty} \dfrac{1}{n} \dfrac{(1-q^n)(1-t^{-n})}{1 + q^n t^{-n}} z^{n} \right)
\end{align}
\smallskip\\
Eq. (\ref{qtVirRelations2}) is the defining relation of an algebraic structure called $q,t$-deformed Virasoro (or just $q$-Virasoro) algebra ${\rm Vir}_{q,t}$ \cite{qtVir}. Abstract definition of this structure requires certain effort because the l.h.s. of (\ref{qtVirRelations2}) is an infinite sum of algebra elements, so it does not make immediate sence to simply define ${\rm Vir}_{q,t}$ as the associative algebra with given generators and relations. Definition can be given either by considering certain completion \cite{qtVirPoisson} or by formalizing extra structures in addition to associative algebra such as the notion of deformed chiral algebra \cite{qtVirDCA}. These details are important but tangential for the purpose of current paper.
\vspace{2ex}

\paragraph{Proposition 6. \cite{qtVir}}In the limit $q = e^{R}, t = e^{\beta R}$ and $R \rightarrow 0$, the operators become scalar, $T_n \rightarrow 2$, and the above relations turn into commutativity relations. The Taylor expansion around $R = 0$ has the form
\begin{align}
T(z) = 2 + \beta \left( z^2 L(z) + \dfrac{(1-\beta)^2}{\beta} \right) R^2 + O(R^4)
\label{qtVirExpansion}
\end{align}
where
\begin{align}
L(z) = \sum\limits_{n \in {\mathbb Z}} \ L_n \ z^{-n}
\end{align}
for some operators $L_n$. Eq. (\ref{qtVirRelations2}) is satisfied for (\ref{qtVirExpansion}) iff operators $L_n$ satisfy the ordinary Virasoro algebra,

\begin{align}
[L_n. L_m] = (n-m) L_{n+m} + c \dfrac{n^3 - n}{12} \delta_{n,-m}
\label{Virasoro2}
\end{align}
\smallskip\\
with central charge $c = 13 - 6 \beta - 6 / \beta$. This justifies the name deformed Virasoro algebra for ${\rm Vir}_{q,t}$.

\newpage

\paragraph{Definition 7.}Let $S(z): {\mathcal F}_{\sigma} \rightarrow {\mathcal F}_{\sigma + 2 \beta}$ be the linear operators defined by

\begin{align}
S(z) = \exp\left( \sum\limits_{k = 1}^{\infty} z^k \dfrac{1 - t^k}{1 - q^k} \dfrac{h_{-k}}{k} \right)
\exp\left( - \sum\limits_{k = 1}^{\infty} z^{-k} (1 + q^k t^{-k}) \dfrac{1 - t^k}{1 - q^k} \dfrac{h_{k}}{k} \right) e^{\beta Q} z^{2 \beta a_0}
\end{align}
\smallskip\\
called screening current operators.
\vspace{2ex}

\paragraph{Proposition 8. \cite{qtVir}}Operators $S(z)$ are distinguished by invariance

\begin{align}
\oint dz \ \big[ T_n , S(z) \big] = 0
\end{align}
\smallskip\\
for closed contours of integration. The contour integrals of $S(z)$ are called screening charge operators.
\vspace{2ex}

\paragraph{Definition 9.}Let $V_{\alpha}(z): {\mathcal F}_{\sigma} \rightarrow {\mathcal F}_{\sigma + \alpha}$ be the linear operators defined by

\begin{align}
V_{\alpha}(z) = \exp\left( \sum\limits_{k = 1}^{\infty} z^k \dfrac{1}{1 + q^k t^{-k}} \dfrac{1 - q^{\alpha k}}{1 - q^k} \dfrac{h_{-k}}{k} \right)
\exp\left( - \sum\limits_{k = 1}^{\infty} z^{-k} \dfrac{1 - q^{-\alpha k}}{1 - q^{-k}} \dfrac{h_{k}}{k} \right) e^{\alpha Q / 2} z^{\alpha a_0}
\end{align}
\smallskip\\
called generic vertex operators.
\vspace{3ex}

\paragraph{Definition 10.}Function

\begin{align}
{\cal B}_{q,t}(x_1, \ldots, x_m) = \oint_{\gamma_1} dz_1 \ldots \oint_{\gamma_N} dz_N \ \big< \alpha_{\infty} \big| \ S(z_1) \ldots S(z_N) \ V_{\alpha_1}(x_1) \ldots V_{\alpha_m}(x_m) \ \big| 0 \big>
\label{BqtDefinition}
\end{align}
\smallskip\\
is called conformal block of ${\rm Vir}_{q,t}$. Here, $\big< \alpha_{\infty} \big| v$ stands for the coefficient of the basis vector $\big| \alpha_{\infty} \big>$ in $v \in {\mathcal F}_{\alpha_{\infty}}$. The function depends on $x_1, \ldots, x_m$ as arguments and $\alpha_1, \ldots \alpha_m, q,t$ as parameters; it also depends on the number and choice of contours $\gamma_1, \ldots, \gamma_N$. Note that by construction $\alpha_{\infty} = \alpha_1 + \ldots + \alpha_m + 2 \beta N$.
\vspace{2ex}

\paragraph{Proposition 11.}In the limit $q = e^{R}, t = e^{\beta R}$ and $R \rightarrow 0$, the function $B_{q,t}(x_1, \ldots, x_m)$ becomes the Dotsenko-Fateev representation of the tree-shaped Virasoro conformal block with $(m+1)$ Verma modules. This is manifest by taking the limit and observing the relation

\begin{align}
V_{\alpha}(z) \ \mathop{\longrightarrow}_{R \rightarrow 0} \ S(z)^{\frac{\alpha}{2\beta}}
\end{align}
\smallskip\\
which is the case for the corresponding Virasoro operators (compare to eqs. (3.1) and (3.33) in \cite{DF}). \linebreak This justifies the name conformal block for the function. Despite no coproduct formalism is available for ${\rm Vir}_{q,t}$, we will sometimes call ${\cal B}_{q,t}(x_1, \ldots, x_m)$ the tree-shaped conformal block, to emphasize the relation to tree-shaped (as opposed to more general graphs with loops) intertwiners of the ordinary Virasoro algebra.

\newpage

\paragraph{Lemma 12.}The following matrix elements can be computed explicitly,

{{\fontsize{9pt}{0pt}{
\begin{align}
\big< S(w) S(z) \big> \equiv \big< 2 \beta \big| S(w) S(z) \big| 0 \big> = w^{2\beta} \ \exp\left( - \sum\limits_{k = 1}^{\infty} (1 + q^{k}t^{-k}) \dfrac{1 - t^k}{1 - q^{k}} \dfrac{z^k w^{-k}}{k} \right) = w^{2\beta} \ \dfrac{\varphi(z/w) \varphi(q/t \ z/w)}{\varphi(q \ z/w) \varphi(t \ z/w)}
\end{align}
\begin{align}
\big< S(z) V_{\alpha}(x) \big> \equiv \big< \alpha + 2 \beta \big| S(z) V_{\alpha}(x) \big| 0 \big> = z^{\alpha} \ \exp\left( - \sum\limits_{k = 1}^{\infty} \dfrac{1 - q^{\alpha k}}{1 - q^{k}} \dfrac{x^k z^{-k}}{k} \right) = z^{\alpha} \ \dfrac{\varphi(x/z)}{\varphi(q^{\alpha} x/z)}
\end{align}
\begin{align}
\big< V_{\gamma}(y) V_{\alpha}(x) \big> \equiv \big< \alpha + \gamma \big| V_{\gamma}(y) V_{\alpha}(x) \big| 0 \big> = y^{\frac{\gamma \alpha}{2 \beta}} \ \exp\left( - \sum\limits_{k = 1}^{\infty} \dfrac{(1-q^{\alpha k})(1-q^{-\gamma k})}{(1-q^{-k})(1-t^k)(1+q^kt^{-k})} \dfrac{x^k y^{-k}}{k} \right)
\end{align}
}}}
\smallskip\\
where $\varphi$ is the Faddeev-Kashaev quantum dilogarithm \cite{FK} closely related to $q$-Gamma and $q$-exponential,

\begin{align}
\varphi(x) = \exp\left( - \sum\limits_{k = 1}^{\infty} \dfrac{1}{1 - q^k} \dfrac{x^k}{k} \right) = \prod\limits_{i = 0}^{\infty} (1 - q^i x) = (x;q)_{\infty} = \dfrac{(q/x)^{\frac{\ln(1-q)}{\ln(q)}} \ (q;q)_{\infty}}{\Gamma_q\left( \frac{\ln(x)}{\ln(q)} \right)} = \dfrac{1}{e_{q}\left( \frac{x}{1-q} \right)}
\label{FKdilog}
\end{align}
\smallskip\\
This is straightforward to compute by using repeatedly the identity $e^{x h} e^{y h^{\dagger}} = e^{xy} e^{y h^{\dagger}} e^{x h}$ for any canonically conjugate operators $[h, h^{\dagger}] = 1$ to normally order the operators featuring in these matrix elements.

\paragraph{Proposition 13.} Conformal block of ${\rm Vir}_{q,t}$ is given by, cf. (\ref{qtConfBlock}):

\begin{align}
{\cal B}_{q,t}(x_1, \ldots, x_m) = \prod\limits_{1 \leq a < b \leq m} \big< V_{\alpha_a}(x_a) V_{\alpha_b}(x_b) \big> \ \cdot \ {\cal I}_{q,t}(x_1, \ldots, x_m)
\label{IntegralRep}
\end{align}
\smallskip\\
where

\begin{align}
{\cal I}_{q,t}(x_1, \ldots, x_m) = \oint_{\gamma_1} dz_1 \ldots \oint_{\gamma_N} dz_N \ \prod\limits_{1 \leq i < j \leq N} \big< S(z_i) S(z_j) \big> \ \prod\limits_{a = 1}^{m} \prod\limits_{i = 1}^{N} \big< S(z_i) V_{\alpha_a}(x_a) \big>
\end{align}
\smallskip\\
This follows from standard results about normal ordering and Wick theorem, expressing the matrix element in the definition of the conformal block as a product of matrix elements with just two operators which we already computed in Lemma 12. The product is taken over all possible choices of two operators, which can be either $S(z_i)$ and $S(z_j)$ (if both operators are chosen from the screening currents) or $S(z_j)$ and $V_{\alpha_a}(x_a)$ (if the first operator is chosen from the vertex operators and the second from screenings) or finally $V_{\alpha_a}(x_a)$ and $V_{\alpha_b}(x_b)$ (if both operators are chosen from the vertex operators). This gives the desired result eq. (\ref{IntegralRep}).

\paragraph{Remark 14.} It is a property of quality special functions to be useful in many different applications. \linebreak
Conformal blocks of ${\rm Vir}_{q,t}$ are equal to instanton partition functions of supersymmetric gauge theories with 8 supercharges in 5 spacetime dimensions, what is known as the generalized AGT conjecture \cite{qtAGT}. \linebreak From the perspective of algebraic geometry, these partition functions represent equivariant integrals \linebreak of appropriate characteristic classes in equivariant $K$-theory of certain algebraic quiver varieties \cite{5dKtheory}. \linebreak Aforementioned equivariant integrals are computed by localization \cite{Localization} resulting in infinite sums over the fixed points of the equivariant torus action. The equality between $K$-theory localization sums and ${\rm Vir}_{q,t}$ conformal blocks is non-trivial and has been checked \cite{qtBlock,qAGT-2,qAGT-3} and proved in different cases \cite{qAGT-Proof-1,qAGT-Proof-2,qAGT-Proof-3,qAGT-Proof-4}.

\newpage

\paragraph{}In this paper, we are interested in a particular special function which represents each of:

\paragraph{$\bullet$} The tree-shaped conformal block of ${\rm Vir}_{q,t}$ algebra (\ref{BqtDefinition}), \\
with five Verma parameters  (i.e. $m = 4$ vertex operators and $\alpha_{\infty}$).

\paragraph{$\bullet$} The instanton partition function of 5-dimensional $SU(2) \times SU(2)$ gauge theory \cite{5dInst}, \\
with 4 fundamental matter multiplets and 1 bifundamental multiplet.

\paragraph{$\bullet$} The equivariant integral in $K$-theory of the ${\cal M}(2) \times {\cal M}(2)$ algebraic variety \cite{5dKtheory}, \\
where {\cal M}(2) is the moduli space of rank 2 torsion free sheaves on ${\mathbb CP}^2$,\\
with integrand given by the Euler class of direct sum of 4 tangent bundles and 1 Ext bundle.

\paragraph{}As explained in Remark 14, from the second and third points of view this function is a localization sum, and equality to the first point of view is not evident. Therefore we give two separate functions, one for the localization sum and one for the conformal block, and connect them momentarily.

\paragraph{Definition 15.} Let ${\cal Z}$ be a function of two complex arguments $\Lambda, x$ given by a formal power series

\begin{align}
{\cal Z}(\Lambda, x)
= \sum\limits_{\nu_1, \nu_2, \mu_1, \mu_2} \ {\mathfrak p}_1^{|\nu_1| + |\nu_2|} \ {\mathfrak p}_2^{|\mu_1| + |\mu_2|}
\ \prod\limits_{1 \leq a,b \leq 2} \dfrac{{\cal N}_{\varnothing,\nu_{b}}\big( v {\mathfrak f}^{+}_a / {\mathfrak n}_b \big)
{\cal N}_{\nu_{a},\mu_{b}}\big( w {\mathfrak n}_a / {\mathfrak m}_b \big)
{\cal N}_{\mu_{b},\varnothing}\big( v {\mathfrak m}_b / {\mathfrak f}^{-}_a \big)
}{ {\cal N}_{\nu_{a},\nu_{b}}\big( {\mathfrak n}_a / {\mathfrak n}_b \big) {\cal N}_{\mu_{a},\mu_{b}}\big( {\mathfrak m}_a / {\mathfrak m}_b \big) }
\end{align}
\smallskip\\
with

\begin{align}
{\cal N}_{\lambda, \eta}\big( z \big) =
\prod\limits_{(i,j) \in \lambda} \left( 1 - z \ q^{\lambda_i - j} t^{\eta^{\prime}_j - i + 1} \right)
\prod\limits_{(i,j) \in \eta} \left( 1 - z \ q^{- \eta_i + j - 1} t^{- \lambda^{\prime}_j + i } \right)
\end{align}
\smallskip\\
Here each of the sums over $\nu_1, \nu_2, \mu_1, \mu_2$ runs through the infinite set $\{ \varnothing, [1], [2], [1,1], [3], [2,1], [1,1,1], \ldots \}$ \linebreak of partitions -- nonincreasing finite sequences of natural numbers. For a given partition $\lambda = [\lambda_1 \geq \lambda_2 \geq \ldots ]$ \linebreak its total size is the sum of its elements, $|\lambda| = \lambda_1 + \lambda_2 + \ldots$ To avoid possible confusion, we note that $\nu_1, \nu_2, \mu_1, \mu_2$ are partitions and their elements are denoted e.g. $(\nu_1)_i$. Each partition has a pictorial \linebreak representation as a Young diagram -- collection of boxes $(i, j) \in {\mathbb Z}_+^2$ such that if box $(i,j)$ is included then $(i-1,j)$ and $(i,j-1)$ are also included in the collection. We use lowercase greek letters to denote partitions and Young diagrams. For a given partition $\lambda$, notation $\prod_{(i,j) \in \lambda}$ stands for the product over boxes of the corresponding Young diagram, i.e. $j$ runs from 1 to $\lambda_i$ and $i$ runs from 1 to the number of elements of $\lambda$. Notation $\lambda^{\prime}$ is the transpose Young diagram or partition, $\lambda^{\prime}_j \equiv $ the number of elements of $\lambda$ which are greater or equal to $j$. Series ${\cal Z}$ is often called Nekrasov function.
\smallskip\\

We choose parametrization

\begin{align}
v = q^{\frac{1}{2}} t^{\frac{-1}{2}}, \ \ \ w = v \phi_1, \ \ \ {\mathfrak p}_1 = v^{-2} T_2 \phi_2 x, \ \ \ {\mathfrak p}_2 = v^{-2} \frac{ T_4 \Lambda}{ \phi_1 x }
\end{align}
\begin{align}
{\mathfrak n}_1 = 1, \ \ {\mathfrak n}_2 = Q, \ \ {\mathfrak m}_1 = 1, \ \ {\mathfrak m}_2 = \phi_1 \phi_2 Q
\end{align}
\begin{align}
{\mathfrak f}^{+}_1 = T_1 Q, \ \ {\mathfrak f}^{+}_2 = T_2^{-1}, \ \ {\mathfrak f}^{-}_1 = T_3^{-1}, \ \ {\mathfrak f}^{-}_2 = T_4 \phi_1 \phi_2 Q
\end{align}
\smallskip\\
so that function ${\cal Z}(\Lambda, x)$ depends on parameters $Q, T_1, T_2, T_3, T_4, \phi_1, \phi_2$ and $q,t$. It is important to emphasize that in this paper ${\cal Z}(\Lambda, x)$ is treated as a formal power series in $\Lambda,x$ where by construction negative powers of $x$ are not excluded but only finitely many can occur for a given finite degree in $\Lambda$. Of course, Nekrasov functions have good and well-studied analytic properties \cite{Convergence} but for now we choose to limit consideration to formal power series, in order to keep presentation as simple and self-contained as possible.

\paragraph{Definition 16.} Let ${\cal J}$ be a function of two complex arguments $\Lambda, x$ given by the integral

\begin{align}
{\cal J}(\Lambda, x) = {\cal I}\big( 0, 1, \Lambda/x, \Lambda \big)
\end{align}
\smallskip\\
with the choice of contours

\begin{align}
{\overrightarrow \gamma} = \big\{ \ \gamma_{[0, \Lambda]}^{N_1}, \ \gamma_{[0, \Lambda/x]}^{N_2}, \ \gamma_{[0, 1]}^{N_3} \ \big\}
\end{align}
\smallskip\\
where $\gamma_{[0,y]}$ denotes a closed contour encircling the line segment $[0,y]$ and $(N_1,N_2,N_3)$ are arbitrary three multiplicities of contours. Function ${\cal J}(\Lambda, x)$ depends on parameters $\alpha_1, \alpha_2, \alpha_3, \alpha_4, N_1, N_2, N_3$ and $q,t$.

\paragraph{Conjecture 17.} Functions ${\cal Z}$ and ${\cal J}$ are related, as formal power series, in the following way

\begin{align}
{\cal J}(\Lambda, x) = \mbox{const} \ \dfrac{\Phi(\phi_1 \phi_2^{-1} q \Lambda x^{-1}) \ \Phi(v^2 T_3 T_4^{-1} q \Lambda x^{-1})}{\Phi(q \Lambda x^{-1})
\ \Phi(v^2 T_3 T_4^{-1} \phi_1 \phi_2^{-1} q \Lambda x^{-1})} \ {\cal Z}\left(\dfrac{\Lambda}{T_2 T_4}, \dfrac{x}{T_2 \phi_2}\right)
\end{align}
\smallskip\\
where we use a convenient notation for the double infinite $q$-Pochhammer product

\begin{align}
\Phi(x) = (x; q, t)_{\infty} = \prod\limits_{i,j = 0}^{\infty} (1 - q^i t^j x) = \exp\left( - \sum\limits_{k = 1}^{\infty} \dfrac{1}{(1 - q^k)(1 - t^k)} \dfrac{x^k}{k} \right)
\end{align}
\smallskip\\
and the map of parameters is given by

\begin{align}
q^{\alpha_1} = v^{-2} T_1 T_2 Q, \ \ \ q^{\alpha_2} = v^{-2}  T_3^{-1} T_4 , \ \ \ q^{\alpha_3} = v^{-2} \phi_2 \phi_1^{-1}, \ \ \ q^{\alpha_4} = v^{-2} T_1 T_2^{-1}
\label{Map1}
\end{align}
\begin{align}
t^{N_1} = v T_1^{-1}, \ \ \ t^{N_2} = v \phi_1, \ \ \ t^{N_3} = v T_3
\label{Map2}
\end{align}
\smallskip\\
The factor "const" may depend on all parameters but not on $\Lambda, x$.

\paragraph{Remark 18.} It is worthwhile to pause again to comment on the meaning of both functions.
First, ${\cal J}(\Lambda, x)$ is the $m = 4$ specialization of (the integral part of) the conformal block for the deformed Virasoro algebra ${\rm Vir}_{q,t}$ with a particular choice of contours parametrized by multiplicities $(N_1,N_2,N_3)$. It is easy to see that this choice of contours is most general in the sense that any other choice can be decomposed in the basis of $(N_1,N_2,N_3)$-contours. The specialization of arguments $(0, 1, \Lambda/x, \Lambda)$ is covenient and not too restrictive.

Second, ${\cal Z}(\Lambda, x)$ is the instanton partition function of 5-dimensional $SU(2) \times SU(2)$ gauge theory with 2 fundamental matter multiplets with mass parameters ${\mathfrak f}^{+}_1, {\mathfrak f}^{+}_2$ charged under the first gauge group and \linebreak 2 multiplets with mass parameters ${\mathfrak f}^{-}_1, {\mathfrak f}^{-}_2$ charged under the second. The gauge groups have Coulomb \linebreak parameters encoded by ${\mathfrak n}_1, {\mathfrak n}_2$ and ${\mathfrak m}_1, {\mathfrak m}_2$ respectively. ${\mathfrak p}_1$ and ${\mathfrak p}_2$ are the instanton parameters of the first and second groups, respectively, and $w$ is up to simple rescaling the mass parameter of the bifundamental.

Third, geometrically, ${\cal Z}(\Lambda, x)$ is an equivariant integral (push-forward) in $K$-theory of the product of two moduli spaces of rank 2 torsion free sheaves on ${\mathbb CP}^2$. The integrand is an appropriate $K$-theory class that represents both bifundamental and fundamental matter in gauge theory. The result of integration is always a sum over fixed points of the torus action, in this case labeled by four partitions $\nu_1, \nu_2, \mu_1, \mu_2$.

\newpage

\section{Non-stationary difference equation}

\paragraph{Definition 19.} Consider the specialization

\begin{align}
\psi(\Lambda,x) = {\cal J}(\Lambda, x) \Big|_{\substack{ N_2 = + 1 \\ \alpha_3 = -\beta}}
\end{align}
\smallskip\\
This is the \textit{degenerate conformal block} with a Verma parameter $\alpha_3 = -\beta$ associated with coordinate $\Lambda / x$ \linebreak and the multiplicity of the corresponding integration contour $\gamma_{[0, \Lambda/x]}$ is specialized to $N_2 = + 1$. \linebreak Note that two conditions are imposed: degeneracy of the Verma parameter and restriction of the contour \linebreak multiplicity (this is well-known at least for ordinary Virasoro, see \cite{DegenFieldReview} for a review). The corresponding \linebreak conditions on the gauge theory side are $\phi_1 = t/v$ and $\phi_2 = v$ which imply relation ${\mathfrak m}_2/{\mathfrak m}_1 = t {\mathfrak n}_2/{\mathfrak n}_1$ between Coulomb parameters and $w = t$ for the bifundamental mass parameter. This specialization of Coulomb parameters and bifundamental mass is known as Higgsing prescription \cite{Higgsing}

\begin{align}
\Psi(\Lambda,x) = {\cal Z}(\Lambda, x) \Big|_{\substack{ \phi_1 = t/v \\ \phi_2 = v}}
\end{align}
\smallskip\\
and the so specialized partition function is the \textit{surface defect wavefunction} in 5-dimensional $SU(2)$ gauge theory with 4 fundamental matter multiplets. The two functions -- degenerate conformal block and surface defect wavefunction -- agree through Conjecture 17, up to the simple prefactor.

\paragraph{Theorem 20.} Function $\Psi(\Lambda,x)$ satisfies, as a formal power series, the non-stationary difference equation:

\begin{align}
\Psi\big(t \Lambda, x\big) = {\cal A}_1(\Lambda, x) \ {\widehat \gamma} \ {\cal A}_2(\Lambda, x) \ {\widehat \gamma} \ {\cal A}_3(\Lambda, x) \ \ \Psi\left(\Lambda, \frac{x}{tqQ}\right)
\label{NonStatEq}
\end{align}
\smallskip\\
where

\begin{align}
{\cal A}_1(\Lambda, x) = \dfrac{1}{\varphi\big(T_1 t v x \big)} \ \dfrac{ \Phi\big( T_3 t^2 v \Lambda x^{-1} \big) }{ \Phi\big( T_3 q v \Lambda x^{-1} \big) } \ \dfrac{ \Phi\big( T_4 t^2 v \Lambda x^{-1} \big) }{ \Phi\big( T_4  t^2 v^{-1} \Lambda x^{-1} \big) }
\end{align}
\begin{align}
{\cal A}_2(\Lambda, x) = \dfrac{ \varphi\big(q T_2 T_3 \Lambda \big) \varphi\big( t T_1 T_4 \Lambda \big) }{ \varphi\big(- T_1 T_2 x \big) \ \varphi\big(- Q^{-1} x \big) \  \varphi\big(- T_3 T_4 Q q t \Lambda x^{-1} \big) \varphi\big(- q \Lambda x^{-1}  \big) }
\end{align}
\begin{align}
{\cal A}_3(\Lambda, x) = \dfrac{1}{\varphi\big(T_2 Q^{-1} q^{-1} v x \big)} \ \dfrac{ \Phi\big( T_3 Q q^2 v \Lambda x^{-1} \big) }{ \Phi\big( T_3 Q q^2 v^{-1} \Lambda x^{-1} \big) } \ \dfrac{ \Phi\big( T_4 Q t^3 v \Lambda x^{-1} \big) }{ \Phi\big( T_4 Q q^2 v^{-1} \Lambda x^{-1} \big) }
\end{align}
\smallskip\\
and the main difference operator ${\widehat \gamma}$ is defined by its action on monomials,

\begin{align}
{\widehat \gamma} \ x^n = q^{\frac{n(n+1)}{2}} \ x^n, \ \ \ n \in {\mathbb Z}
\end{align}
\smallskip\\
The simplest operator with this property is ${\widehat \gamma} = q^{\partial(\partial+1)/2}$ where $\partial = \partial_{\ln x} = x \partial_x $. The reason why this difference equation is called non-stationary is because the l.h.s. is shifted with respect to $\Lambda$ and the magnitude of the shift is an independent variable $t$, therefore one may think of this as a non-stationary Schrodinger equation describing discrete time-evolution of a quantum-mechanical wavefunction.

\paragraph{}We checked eq.(\ref{NonStatEq}) in power series expansion up to degree $O\big(\Lambda^6, x^{10}\big)$ for all parameters generic. \linebreak Before proving Theorem 20 below, we consider first an interesting limit when parameters $T_1,T_2,T_3,T_4 \rightarrow 0$.

\pagebreak

\paragraph{Corollary 21.} Theorem 20 implies that function

\begin{align}
\Psi_{\rm Toda}(\Lambda,x) \ = \ \lim\limits_{\substack{ { T_1, T_2, T_3, T_4 \rightarrow 0 } \\ {\tiny{ T_1/T_2, T_3/T_4 \ {\rm fixed}}} }} \ \Psi(\Lambda,x)
\end{align}
\smallskip\\
satisfies a much simpler equation,

{\fontsize{9pt}{0pt}
\begin{align}
\Psi_{\rm Toda}\big(t \Lambda, x\big) \ = \ {\hat H} \Psi_{\rm Toda}\big(\Lambda, x\big), \ \ \ \ \ {\hat H} f(\Lambda, x) \ \equiv \ {\widehat \gamma} \ \prod\limits_{n = 0}^{\infty} \left( 1 + q^n \frac{x}{Q} \right)^{-1} \left( 1 + q^{n+1} \frac{\Lambda}{x} \right)^{-1}  \ {\widehat \gamma} \ \ f\left( \Lambda, \frac{x}{tqQ}\right)
\label{NonStatEqToda}
\end{align}}
\smallskip\\
Eq. (\ref{NonStatEqToda}) is identical to the non-stationary relativistic (affine) Toda equation considered in \cite{Shiraishi1}, eq. (11). \linebreak This choice of name for the equation is supported by the following

\paragraph{Proposition 22.}Operator ${\hat H}$ and the quantum Hamiltonian of the 2-body relativistic affine Toda system

\begin{align}
{\hat H}_{\rm Toda} f(\Lambda, x) \ \equiv \ f(\Lambda, q x) + t Q f(\Lambda, q^{-1} x) + t x f(\Lambda, x) + \Lambda x^{-1} f(\Lambda, x)
\label{TodaHamiltonian}
\end{align}
\smallskip\\
commute with each other:
\begin{align}
[ {\hat H}, {\hat H}_{\rm Toda} ] = 0
\end{align}
\paragraph{Proof.}This statement appears, for example, as Proposition 6.2. in \cite{Shiraishi1}, however here we choose to present it in some detail for completeness. We rely upon commutation relation between ${\widehat \gamma}$ and basic operators

\begin{align}
{\hat p} f(\Lambda, x) = f(\Lambda, q x), \ \ \ {\hat x} f(\Lambda, x) = x f(\Lambda, q x), \ \ \ {\hat p} {\hat x} = q {\hat x} {\hat p}
\end{align}
\smallskip\\
These operators satisfy ${\hat p} {\widehat \gamma} = {\widehat \gamma} {\hat p}$ and ${\widehat \gamma} {\hat x} = {\hat p} {\hat x} {\widehat \gamma}$. Using these relations, we find

{\fontsize{9pt}{0pt}{
\begin{align}
{\hat H}_{\rm Toda} {\hat H} \ f(\Lambda,x) = {\widehat \gamma} \ \Big( {\hat p} + t Q {\hat p}^{-1} + t {\hat p}^{-1} {\hat x} + \Lambda {\hat x}^{-1} {\hat p} \Big) \ \varphi( -x/Q )^{-1} \varphi( -q\Lambda/x )^{-1} \ {\widehat \gamma} \ \ f\left( \Lambda, \frac{x}{tqQ}\right)
\end{align}
\begin{align}
{\hat H} {\hat H}_{\rm Toda} \ f(\Lambda,x) = {\widehat \gamma} \ \varphi( -x/Q )^{-1} \varphi( -q\Lambda/x )^{-1}  \
\Big( {\hat p} + t Q {\hat p}^{-1} + Q^{-1} q^{-1} {\hat p} {\hat x} + q t Q \Lambda {\hat x}^{-1} {\hat p}^{-1} \Big) \
{\widehat \gamma} \ \ f\left( \Lambda, \frac{x}{tqQ}\right)
\end{align}
}}\smallskip\\
One can see that $[ {\hat H}, {\hat H}_{\rm Toda} ] \ f(\Lambda,x) = {\widehat \gamma} \ {\widehat O} \ {\widehat \gamma} \ \ f\left( \Lambda, \frac{x}{tqQ}\right)$ with
\begin{align}
{\widehat O} \ = \ & \Big[ (1 + \Lambda x^{-1}) {\hat p} + ( t Q + t q^{-1} x ) {\hat p}^{-1} \Big] \ \varphi( -x/Q )^{-1} \varphi( -q\Lambda/x )^{-1} - \emph{} \\
& \emph{} - \varphi( -x/Q )^{-1} \varphi( -q\Lambda/x )^{-1}
\Big[ ( 1 + x / Q ) {\hat p} + ( t Q + q t Q \Lambda x^{-1} ) {\hat p}^{-1} \Big]
\end{align}
Using property $\varphi(q x) = (1-x)^{-1} \varphi(x)$ we find
\begin{align}
{\widehat O} \ = \ & \varphi( -x/Q )^{-1} \varphi( -q\Lambda/x )^{-1} \left[ (1 + \Lambda x^{-1}) \dfrac{1 + x/Q}{1 + \Lambda / x} - ( 1 + x/Q) \right] {\hat p} + \emph{} \\
& \emph{} +
\varphi( -x/Q )^{-1} \varphi( -q\Lambda/x )^{-1} \left[ (t Q + t q^{-1} x) \dfrac{1 + q \Lambda x^{-1}}{1 + q^{-1} x/Q} - ( t Q + t q Q \Lambda x^{-1} ) \right] {\hat p}^{-1} = 0
\end{align}
where each of the rational functions in brackets vanishes, so ${\widehat O}$ is a zero operator, completing the proof. \qedsymbol

\paragraph{Remark 23.} We contend it is highly likely that the operator in the r.h.s. of eq. (\ref{NonStatEq}) commutes with the appropriate deformation of Toda Hamiltonian (\ref{TodaHamiltonian}). From general principles \cite{SW-Integ-Mir} this deformation should represent the Hamiltonian of XXZ integrable model. We leave this as a straightforward exercise.

\pagebreak

\paragraph{Proof of Theorem 20.}In this paper we only prove Theorem 20 for a concrete choice
\begin{align}
(\alpha_2,\alpha_3,\alpha_4; N_1,N_2,N_3) = \big(-\beta, -\beta, -\beta; 1,1,0 \big)
\label{ParameterChoice}
\end{align}
equivalently
\begin{align}
(T_1, T_2, T_3, T_4, \phi_1, \phi_2) = \big( v t^{-1}, \ v^{-1}, \ v^{-1}, \ v t^{-1}, \ t v^{-1}, \ v \big)
\label{ParameterChoice2}
\end{align}
which is distinguished by all three parameters $\alpha_1, \alpha_2, \alpha_3$ being degenerate. Note that $\alpha_1$ is still generic. \linebreak This choice is intended as a simple starting point for the upcoming general proof, that we hope will be found in near future. We argue below that already this case is not entirely trivial and produces interesting hints.

Generally $\Psi(\Lambda, x)$ is a sum over 4 partitions $\nu_1,\nu_2,\mu_1,\mu_2$, however at the point (\ref{ParameterChoice2}) its summand \linebreak contains a particular product ${\cal N}_{\varnothing, \nu_1}\big( q/t \big) {\cal N}_{\varnothing, \nu_2}\big( q/t^2 \big) {\cal N}_{\mu_1, \varnothing}\big( 1 \big) {\cal N}_{\mu_2, \varnothing}\big( t \big)$ what leads to drastic simplification, because this product is only non-vanishing for $(\nu_1,\nu_2,\mu_1,\mu_2) = (\varnothing,[n],\varnothing,[m])$. Therefore the sum over 4 partitions reduces, in this case, to a sum over two non-negative integers $n$ and $m$. Short calculation gives

\begin{align}
\Psi(\Lambda, x) \Big|_{{\rm at} \ (\ref{ParameterChoice2})} \ \equiv \ U(\Lambda, x) \ = \ \sum\limits_{n,m = 0}^{\infty} \ q^{n m + n} t^{-n} \ \dfrac{(t Q)_n (qQ/t)_m}{(q Q)_m ( qQ/t )_n} \ (q^{-m})_n (v^2 q^{-n})_m \ \dfrac{x^{n} }{(q)_n} \ \dfrac{(\Lambda/x)^{m}}{(q)_m}
\label{UofX}
\end{align}
\smallskip\\
where we use a standard notation for $q$-Pochhammer symbols $(z)_n = \prod_{i = 0}^{n-1}(1 - q^i z)$. Function $U(\Lambda,x)$ is in fact identical to one of the most famous special functions -- the Macdonald polynomial; indeed, one has

\begin{align}
\nonumber \dfrac{t(1-xt)(1-\Lambda)}{(1-x)(t-\Lambda)} \ U\big( q \Lambda, q x \big) +
\dfrac{(q \Lambda - t^2 x)(x - t)}{t(1-x)(\Lambda - t x)(x - 1)} \ U\big( \Lambda, x/q \big) \\ +
\dfrac{t(\Lambda - q x)(\Lambda - t^2)}{q Q (\Lambda - t x)(\Lambda - t)} \ U\big( \Lambda/q, x \big) = (1 + t + t/Q) \ U\big( \Lambda, x \big)
\label{MacdonaldOperator}
\end{align}
\smallskip\\
which is completely straightforward to prove from the series definition eq. ({\ref{UofX}}) by rewriting ({\ref{MacdonaldOperator}}) as a rational recursive relation for coefficients. It is not hard to recognize in eq. (\ref{MacdonaldOperator}) the Macdonald \cite{Macd} or Ruijsenaars-Schneider \cite{RS} difference equation for Macdonald polynomials of $U(3)$, i.e. in 3 variables:

\begin{align}
U\big(\Lambda, x\big) \Big|_{Q = q^{-r} t^{-1}} = \dfrac{\varphi\big(q t^{-2} \Lambda x^{-1} \big)}{\varphi\big(\Lambda x^{-1} \big)} \ P_{[r, 0, 0]}\big( 1, \ t^{-1} \Lambda, \ t^{-1} \Lambda/x \big)
\label{MacdonaldSolution}
\end{align}
\smallskip\\
where $P$ are the standard Macdonald polynomials labeled by symmetric representations $[r,0,0]$ of $U(3)$:

\begin{align}
P_{[1, 0, 0]}(x_1,x_2,x_3) \ = \ & x_1 + x_2 + x_3 \\ \nonumber
\\
P_{[2, 0, 0]}(x_1,x_2,x_3) \ = \ & x_1^2 +  x_2^2 + x_3^2 + \dfrac{(1+q)(1-t)}{1-qt} (x_1 x_2 + x_1 x_3 + x_2 x_3) \\ \nonumber
\\
\nonumber P_{[3, 0, 0]}(x_1,x_2,x_3) \ = \ & x_1^3 +  x_2^3 + x_3^3 + \dfrac{(1+q+q^2)(1-t)}{1-q^2t} (x^2_1 x_2 + x_1 x_2^2 + x_1^2 x_3 + \emph{} \\
& \emph{} + x_1 x_3^2 + x_2^2 x_3 + x_2 x_3^2) + \dfrac{(1+q)(1+q+q^2)(1-t)^2}{(1-qt)(1-q^2t)} x_1 x_2 x_3, \ \ \ \ \ \ldots
\label{MacdonaldPolynomials3}
\end{align}
\smallskip\\
and so forth. This clarifies the nature of the function $\Psi(\Lambda, x)$ in this case.

\newpage

As for the equation that it satisfies, it is important to note that equation (\ref{MacdonaldOperator}) is the most basic, but not the one best suitable for generalization to bigger integers $N_1,N_2,N_3$. Instead, we want to use this very simple example to address equation (\ref{NonStatEq}) which is more complicated, but instead fully general. We have

\begin{align}
{\cal A}_1(\Lambda, x) \Big|_{{\rm at} \ (\ref{ParameterChoice2})} \equiv u_1(\Lambda,x) = \dfrac{1}{\varphi\big(q t^{-1} x\big)\varphi\big(t \Lambda x^{-1} \big)}
\end{align}
\begin{align}
{\cal A}_2(\Lambda, x) \Big|_{{\rm at} \ (\ref{ParameterChoice2})} \equiv u_2(\Lambda,x) = \dfrac{ \varphi\big(t \Lambda \big) \varphi\big( q t^{-2} \Lambda \big) }{ \varphi\big(- t^{-1} x \big) \ \varphi\big(- Q^{-1} x \big) \  \varphi\big(- q Q \Lambda x^{-1} \big) \varphi\big(- q \Lambda x^{-1}  \big) }
\end{align}
\begin{align}
{\cal A}_3(\Lambda, x) \Big|_{{\rm at} \ (\ref{ParameterChoice2})} \equiv u_3(\Lambda,x) = \dfrac{1}{\varphi\big(q^{-1} Q^{-1} x\big)\varphi\big(q^2 t^{-1} Q \Lambda x^{-1} \big)}
\end{align}
\smallskip\\
and therefore equation (\ref{NonStatEq}) takes form

\begin{align}
\label{EqnV1} V(\Lambda, x) \ \equiv \ \dfrac{ \varphi\big(- t^{-1} x \big) \  \varphi\big(- q \Lambda x^{-1} \big) }{ \varphi\big( q t^{-2} \Lambda \big) } \ \
{\widehat \gamma}^{-1} \ \varphi\big(q t^{-1} x\big) \ \varphi\big(t \Lambda x^{-1} \big) \ \ U\big(t \Lambda, x\big) = \emph{} \\
\nonumber \\
\label{EqnV2} \emph{} = \dfrac{ \varphi\big( t\Lambda \big) }{ \varphi\big(- Q^{-1} x \big) \varphi\big(- q Q \Lambda x^{-1} \big) } \ \
{\widehat \gamma}^{+1} \ \dfrac{ 1 }{ \varphi\big(q^{-1} Q^{-1} x\big) \varphi\big(q^2 t^{-1} Q \Lambda x^{-1} \big) } \ U\left( \Lambda, \frac{x}{tqQ}\right)
\end{align}
\smallskip\\
\textit{This is what we want to prove.} Our strategy is based on proving an explicit formula for $V(\Lambda,x)$:

\begin{align}
V(\Lambda, x) \ = \ \sum\limits_{n,m = 0}^{\infty} \ q^{m} (-Q)^{m-n} \ \dfrac{(t Q)_n (t)_m }{(q Q)_m (t)_n} \ \dfrac{ (q^{-m})_n (v^2 q^{-n})_m }{(v^2 q^{-m})_n (t^2 q^{-n})_m} (q t^{-2})_n (t^2)_m \ \dfrac{x^{n} }{(q)_n} \ \dfrac{(\Lambda/x)^{m}}{(q)_m}
\label{EqnV3}
\end{align}
\smallskip\\
which satisfies both eq. (\ref{EqnV1}) and (\ref{EqnV2}). We find it quite remarkable that such explicit formula for $V(\Lambda,x)$ exists at all, allowing to speculate that non-stationary difference equation admits a "square root" (can be split into two equations) and we comment on this in more detail in the Conclusion. Here, let us proceed to prove the theorem, starting from the first of two equations:

\begin{align}
{\widehat \gamma} \ \varphi\big( q t^{-2} \Lambda \big) \ \varphi\big(- t^{-1} x \big)^{-1} \ \varphi\big(- q \Lambda x^{-1} \big)^{-1} \ V(\Lambda, x) \ = \ \varphi\big(q t^{-1} x\big) \ \varphi\big(t \Lambda x^{-1} \big) \ \ U\big(t \Lambda, x\big)
\label{ProofEquation1}
\end{align}
\smallskip\\
Both l.h.s. and r.h.s. are formal power series in $\Lambda,x$ where each given term is a rational function of $Q$ whose degree in $Q$ is bounded by the degree in $\Lambda,x$. Two such series coincide for generic $Q$ if they coincide for all

\begin{align}
Q = q^{-r} t^{-1}, \ \ \ r = 0,1,2,3,\ldots
\end{align}
\smallskip\\
This allows to prove the desired equation using generating functions. For r.h.s. of eq. (\ref{ProofEquation1}) it takes form

{\fontsize{9pt}{0pt}
\begin{align}
\sum\limits_{r=0}^{\infty} \ z^r \ \dfrac{(t)_r}{(q)_r} \ \varphi\big(q t^{-1} x\big) \ \varphi\big(t \Lambda x^{-1} \big) U\big(t \Lambda, x\big) \Big|_{Q = q^{-r} t^{-1}} = \varphi\big(q t^{-1} x\big) \ \varphi\big(q t^{-1} \Lambda x^{-1} \big) \ \dfrac{\varphi\big( t z \big) \varphi\big( t \Lambda z \big) \varphi\big( t \Lambda x^{-1} z\big)}{\varphi\big( z\big) \varphi\big( \Lambda z\big) \varphi\big( \Lambda x^{-1} z\big)}
\label{GenFuncLHS1}
\end{align}}
\smallskip\\
as a corollary of the well-known generating function for Macdonald polynomials. Generating function for the l.h.s. of \pagebreak eq. (\ref{ProofEquation1}) is from the first glance somewhat more involved, but we can still compute it using the basic $q$-binomial theorem

\begin{align}
\dfrac{\varphi\big(x\big)}{\varphi\big(x/y\big)} = \sum\limits_{n = 0}^{\infty} \ \dfrac{(y)_n}{(q)_n} \ (x/y)^n
\end{align}
\smallskip\\
to perform the summation over $r$:

{\fontsize{9pt}{0pt}{
\begin{align}
\nonumber & \varphi\big( q t^{-2} \Lambda \big) \ {\widehat \gamma} \ \ \sum\limits_{r=0}^{\infty} \ z^r \ \dfrac{(t)_r}{(q)_r} \ \varphi\big(- t^{-1} x \big)^{-1} \ \varphi\big(- q \Lambda x^{-1} \big)^{-1} \ V(\Lambda, x) \Big|_{Q = q^{-r} t^{-1}} = \emph{} \\
\nonumber & = \varphi\big( q t^{-2} \Lambda \big) \ \dfrac{\varphi\big(t z\big)}{\varphi\big(z\big)} \ {\widehat \gamma} \varphi\big(- t^{-1} x \big)^{-1} \ \varphi\big(- q \Lambda x^{-1} \big)^{-1}
\sum\limits_{n,m = 0}^{\infty} \ z^{n+m} \ t^n q^{n(n-1)/2 - m(m-1)/2} \times \emph{} \\
& \emph{} \times \dfrac{ (q^{-m})_n (v^2 z^{-1} q^{-n})_m }{(v^2 q^{-m})_n (t^2 q^{-n})_m} \dfrac{(t)_m}{(t z)_n} (q t^{-2})_n (t^2)_m \ \dfrac{x^{n} }{(q)_n} \ \dfrac{(\Lambda/x)^{m}}{(q)_m} \ \mathop{=}^{\star} \ \emph{}
\label{Unwieldy}
\end{align}}}
\smallskip\\
This unwieldy expression will momentarily be simplified, once we apply the operator ${\widehat \gamma}$ explicitly via

\begin{align}
{\widehat \gamma} \varphi\big(- t^{-1} x \big)^{-1} \ \varphi\big(- q \Lambda x^{-1} \big)^{-1} \ x^{n-m} = \varphi\big( q t^{-1} \Lambda \big)^{-1} \varphi\big( q^{1 + n - m} t^{-1} x \big)\varphi\big( q^{1 + m - n} \Lambda x^{-1} \big) \ {\widehat \gamma} \ x^{n-m}
\label{FormulaGamma}
\end{align}
\smallskip\\
which is a direct corollary of the definition of ${\widehat \gamma}$ and elementary series expansions

\begin{align}
\varphi\big(x\big) = \sum\limits_{n = 0}^{\infty} \ q^{n(n+1)/2} \ \dfrac{(-x)^n}{(q)_n}, \ \ \ \varphi\big(x\big)^{-1} = \sum\limits_{n = 0}^{\infty} \dfrac{x^n}{(q)_n}
\end{align}
\smallskip\\
The application of (\ref{FormulaGamma}) to (\ref{Unwieldy}) leads to a major simplification,

\begin{align}
\nonumber \emph{} \ \mathop{=}^{\star} \ & \dfrac{\varphi\big( q t^{-2} \Lambda \big)}{\varphi\big( q t^{-1} \Lambda \big)} \ \dfrac{\varphi\big(t z\big)}{\varphi\big(z\big)} \ \sum\limits_{n,m = 0}^{\infty} \ z^{n+m} \ t^n q^{n^2 - nm } \dfrac{ (q^{-m})_n (v^2 z^{-1} q^{-n})_m }{(v^2 q^{-m})_n (t^2 q^{-n})_m} \times \emph{} \\
& \emph{} \times \dfrac{(t)_m}{(t z)_n} (q t^{-2})_n (t^2)_m \ \varphi\big( q^{1 + n - m} t^{-1} x \big)\varphi\big( q^{1 + m - n} \Lambda x^{-1} \big) \dfrac{x^{n} }{(q)_n} \ \dfrac{(\Lambda/x)^{m}}{(q)_m}
\end{align}
\smallskip\\
Equating this to (\ref{GenFuncLHS1}), we get the relation between generating functions:

{\fontsize{7pt}{0pt}{
\begin{align}
\nonumber & \sum\limits_{n,m = 0}^{\infty} \ z^{n+m} \ t^n q^{n^2 - n m } \dfrac{ (q^{-m})_n (v^2 z^{-1} q^{-n})_m }{(v^2 q^{-m})_n (t^2 q^{-n})_m} \times \emph{} \\
& \emph{} \times \dfrac{(t)_m}{(t z)_n} (q t^{-2})_n (t^2)_m \ \dfrac{\varphi\big( q^{1 + n - m} t^{-1} x \big)\varphi\big( q^{1 + m - n} \Lambda x^{-1} \big)}{\varphi\big(q t^{-1} x\big) \ \varphi\big(q t^{-1} \Lambda x^{-1} \big) } \dfrac{x^{n} }{(q)_n} \ \dfrac{(\Lambda/x)^{m}}{(q)_m} = \dfrac{\varphi\big( q t^{-1} \Lambda \big)}{\varphi\big( q t^{-2} \Lambda \big)} \ \dfrac{\varphi\big( t \Lambda z \big) \varphi\big( t \Lambda x^{-1} z\big)}{\varphi\big( \Lambda z\big) \varphi\big( \Lambda x^{-1} z\big)}
\end{align}}}
\smallskip\\
Changing summation index from $n$ to $k = m - n$, and simplifying the Pochhammer symbols, we find

{\fontsize{7pt}{0pt}{\begin{align}
\sum\limits_{k = 0}^{\infty} q^{k^2} \left(\frac{t}{x}\right)^k \ \ \dfrac{\varphi\big( q^{1 - k} t^{-1} x \big)\varphi\big( q^{1 + k} \Lambda x^{-1} \big)}{\varphi\big(q t^{-1} x\big) \ \varphi\big(q t^{-1} \Lambda x^{-1} \big) } \ \dfrac{(q t^{-1} z^{-1})_{k} (t)_k }{ (q)_k (t^2)_k }
\sum\limits_{m = k}^{\infty} \left(\frac{z \Lambda}{t q^k}\right)^m \ \dfrac{(t^2)_m}{(q)_{m-k}}
 = \dfrac{\varphi\big( q t^{-1} \Lambda \big)}{\varphi\big( q t^{-2} \Lambda \big)} \ \dfrac{\varphi\big( t \Lambda z \big) \varphi\big( t \Lambda x^{-1} z\big)}{\varphi\big( \Lambda z\big) \varphi\big( \Lambda x^{-1} z\big)}
\end{align}}}
\smallskip\\
Using $q$-binomial theorem again, the sum over $m$ can be finally removed and \pagebreak all $\varphi$-functions under summation can be traded for Pochhammers, resulting in the ultimate form of the equation that we need to prove:

\begin{align}
\sum\limits_{k = 0}^{\infty} \ q^k \ \dfrac{(t)_k (t x^{-1})_k (q t^{-1} z^{-1})_k}{ (q t z^{-1} \Lambda^{-1} )_k (q \Lambda x^{-1})_k (q)_k } \ = \ \dfrac{\varphi\big( q t^{-1} \Lambda x^{-1} \big)  \varphi\big( t^{-1} \Lambda z \big) \varphi\big( t \Lambda x^{-1} z\big) \varphi\big( q t^{-1} \Lambda \big)}{\varphi\big( q \Lambda x^{-1} \big) \varphi\big( \Lambda z\big) \varphi\big( \Lambda x^{-1} z\big) \varphi\big( q t^{-2} \Lambda \big)}
\label{qSeriesIdentity}
\end{align}
\smallskip\\
At this point we are essentially done, because this is a well-known identity in the theory of $q$-deformed (basic) hypergeometric functions in one variable \cite{Gasper-Rahman}. Recall that all equations in this section are understood as equalities between formal power series in $\Lambda, x$ such that negative powers of $x$ are permitted in the expansion, but for any given term the negative degree in $x$ cannot exceed the positive degree in $\Lambda$. When generating functions in $z$ are considered, each term in the $\Lambda, x$ expansion is understood as a formal power series in positive powers of $z$. For example, eq. (\ref{qSeriesIdentity}) is the following expansion

\begin{align}
1 + \dfrac{(1-t)(t-x)(tz-q)}{x t^2 (1-q)} \Lambda + O(\Lambda^2) \ = \ \dfrac{\varphi\big( q t^{-1} \Lambda x^{-1} \big)  \varphi\big( t^{-1} \Lambda z \big) \varphi\big( t \Lambda x^{-1} z\big) \varphi\big( q t^{-1} \Lambda \big)}{\varphi\big( q \Lambda x^{-1} \big) \varphi\big( \Lambda z\big) \varphi\big( \Lambda x^{-1} z\big) \varphi\big( q t^{-2} \Lambda \big)}
\end{align}
\smallskip\\
By construction, $k$-th term in the $\Lambda$ expansion is a Laurent polynomial in $t$ of order $-2k$ and degree $+k$. Therefore, as a power series equality, eq.(\ref{qSeriesIdentity}) holds for generic $t$ if it holds for all

\begin{align}
t = q^{-n}, \ \ \ n = 1,2,3,\ldots
\end{align}
\smallskip\\
At these points eq.(\ref{qSeriesIdentity}) becomes

\begin{align}
\sum\limits_{k = 0}^{\infty} \ q^k \ \dfrac{( q^{-n} )_k ({\mathit a})_k ({\mathit b})_k}{ ({\mathit c})_k (d)_k (q)_k } \ = \ \dfrac{({\mathit c}/{\mathit a})_n ({\mathit c}/{\mathit b})_n}{({\mathit c})_n ({\mathit c} {\mathit a}^{-1} {\mathit b}^{-1} )_n}
\label{qSaalschutz}
\end{align}
\smallskip\\
where, we remind, the $q$-Pochhammer symbols are $(z)_n = \prod_{i = 0}^{n-1}(1 - q^i z)$, and parameters

\begin{align}
{\mathit a} = t x^{-1}, \ \ \ {\mathit b} = q t^{-1} z^{-1}, \ \ \ {\mathit c} = q t z^{-1} \Lambda^{-1}, \ \ \ {\mathit d} = q \Lambda x^{-1}
\end{align}
\smallskip\\
satisfy the balancing condition $ q t {\mathit a} {\mathit b} = {\mathit c} {\mathit d}$. Note that it is not any problem for this argument that $t = q^{-n}$ is a pole of many important functions, including Macdonald polynomials (\ref{MacdonaldPolynomials3}) we used just above, and was explicitly prohibited when we constructed conformal blocks from bosonic representations in Section 2. \linebreak The specialization $t = q^{-n}$ is only used here as a technical tool to prove the isolated identity (\ref{qSeriesIdentity}) and for this purpose it is perfectly fine. This completes the (first half of the) proof, since identity (\ref{qSaalschutz}) is the well-established $q$-Saalsch\"utz identity \cite{Gasper-Rahman} due to Jackson \cite{qSaalschutzJackson}. For the second half of the proof, we would need to also prove the second of two main difference equations, that is eq. (\ref{EqnV2}),

\begin{align}
\varphi\big( t\Lambda \big) \ {\widehat \gamma}^{-1} \ \varphi\big(- Q^{-1} x \big) \varphi\big(- q Q \Lambda x^{-1} \big) \ V(\Lambda, x) \ = \ \dfrac{ 1 }{ \varphi\big(q^{-1} Q^{-1} x\big) \varphi\big(q^2 t^{-1} Q \Lambda x^{-1} \big) } \ U\left( \Lambda, \frac{x}{tqQ}\right)
\label{ProofEquation2}
\end{align}
\smallskip\\
which we skip, because it is a step-by-step repetition of the proof for (\ref{ProofEquation1}) presented in some detail above. \qedsymbol
\smallskip\\

Before proceeding to Conclusion, we comment on possibility to solve non-stationary difference equations. We consider only the case $T_1,T_2,T_3,T_4 \rightarrow 0$ for brevity, generalization is expected to be immediate.

\newpage

\paragraph{Theorem 24.}Function $\Psi_{\rm Toda}\big(\Lambda, x\big)$ admits the following representation

{\fontsize{8pt}{0pt}
\begin{align}
\Psi_{\rm Toda}\big(\Lambda, x\big) \ = \ \prod\limits_{m = 0}^{\infty} \ \left\{ \ {\widehat \gamma} \ \prod\limits_{n = 0}^{\infty} \left( 1 + q^{n-m} t^{-m} \frac{x}{Q^{m+1}} \right)^{-1} \left( 1 + q^{n+m+1} t^{-1} \frac{\Lambda Q^{m}}{x} \right)^{-1}  \ {\widehat \gamma} \ \ \right\} \ \cdot \ 1
\label{Wrepresentation}
\end{align}}
\smallskip\\
as the action of explicit operator on an elementary function, sometimes called cut-and-join representation \cite{CutAndJoin1,CutAndJoin2,CutAndJoin2} or ${\cal W}$-representation \cite{Wrep1,Wrep2,Wrep3}. This can be regarded as a rather direct solution of eq. (\ref{NonStatEqToda}).

\paragraph{Proof.}Provided Theorem 20, function $\Psi_{\rm Toda}\big(\Lambda, x\big)$ satisfies the equation (\ref{NonStatEqToda}), which can be written as

{\fontsize{8pt}{0pt}
\begin{align}
\Psi_{\rm Toda}\big(\Lambda, x\big) \ = \ {\widehat \gamma} \ \prod\limits_{n = 0}^{\infty} \left( 1 + q^n \frac{x}{Q} \right)^{-1} \left( 1 + q^{n+1} t^{-1} \frac{\Lambda}{x} \right)^{-1}  \ {\widehat \gamma} \ \ \Psi_{\rm Toda}\left( \frac{\Lambda}{t}, \frac{x}{tqQ}\right)
\end{align}}
\smallskip\\
Iterating the same formula by applying to the $\Psi$-function in the r.h.s. $M$ times, we find

{\fontsize{8pt}{0pt}
\begin{align}
\hspace{-3ex} \ \Psi_{\rm Toda}\big(\Lambda, x\big) \ = \ \prod\limits_{m = 0}^{M} \left\{ \ {\widehat \gamma} \prod\limits_{n = 0}^{\infty} \left( 1 + q^{n-m} t^{-m} \frac{x}{Q^{m+1}} \right)^{-1} \left( 1 + q^{n+m+1} t^{-1} \frac{\Lambda Q^{m}}{x} \right)^{-1}  \ {\widehat \gamma} \ \right\} \ \cdot \ \Psi_{\rm Toda}\left( \frac{\Lambda}{t^{M+1}}, \frac{x}{(tqQ)^{M+1}}\right)
\end{align}}
\smallskip\\
At least for parameters $q,t,Q$ in the region

\begin{align}
t \in (1, \infty) \in {\mathbb R}, \ \ \ q t Q \in (1, \infty) \in {\mathbb R}, \ \ \ q^{-1} Q^{-1} \in (1, \infty) \in {\mathbb R}
\end{align}
\smallskip\\
the function in the r.h.s. satisfies, as a formal power series,

\begin{align}
\lim\limits_{M \rightarrow \infty} \ \Psi_{\rm Toda}\left( \frac{\Lambda}{t^{M+1}}, \frac{x}{(tqQ)^{M+1}}\right) = 1
\end{align}
\smallskip\\
and similarly the $M \rightarrow \infty$ limit of the product in the r.h.s. is an infinite product that is well-defined as an operator on formal power series in $\Lambda,x$. Taking the limit $M \rightarrow \infty$ gives eq. (\ref{Wrepresentation}). \qedsymbol

\section{Conclusion}

Non-stationary difference equations associated with many-body quantum integrable systems are interesting from many perspectives. From pure integrability or even more generally quantum mechanics standpoint, it is remarkable that right hand sides of these equations are \textit{defined by infinite products}, operators of the form
\begin{align}
{\widehat \gamma} \ \prod\limits_{n = 0}^{\infty} \left( 1 + q^n x \right)^{-1} \left( 1 + q^{n+1} \Lambda x^{-1} \right)^{-1}  \ {\widehat \gamma}
\end{align}
which commute with the Hamiltonian operators that define stationary difference equations, such as
\begin{align}
{\widehat p} + {\widehat p}^{-1} + x + \Lambda x^{-1}
\end{align}
where both examples represent affine relativistic Toda with 2 particles (1 center-of-mass degree of freedom). \linebreak Of course there is nothing surprising in the fact that some operators commute with the Hamiltonians, because one can construct polynomials and more general functions of Hamiltonians, but what is surprising is that \linebreak particular function of the Hamiltonians takes such an explicit product form. To the best of our knowledge, \linebreak one of the first occurrences of this kind of operators was eq. (3.24) in \cite{KashaevOperators}. We believe this should hold generally, well beyond the Toda case, meaning for every integrable system there exist both standard \linebreak stationary Hamiltonians and their "big" non-stationary counterparts, that commute with each other, but have drastically different properties. In some aspects, the latter operators can be simpler than the original ones: for example, we showed in the proof of Theorem 20 that non-stationary difference equation admits a \textit{"square root" type decomposition} into two parts (eqs. (\ref{EqnV1}) and (\ref{EqnV2}), provided (\ref{EqnV3})) while we are unaware of the parallel phenomenon for the regular stationary Toda Hamiltonian.

This touches a plethora of interrelated subjects: for one, such decomposition is reminiscent of modular $SL(2,{\mathbb Z})$ identities of the form $S^2 = \mbox{id}$ for the Fourier transform operators and $(ST)^3 = \mbox{id}$ for Dehn twists for Macdonald polynomials at roots of unity \cite{inner} further generalized to \textit{modular $SL(3,{\mathbb Z})$ identities} \cite{SL3Z} for affine and elliptic Macdonald functions. These resemblances are unlikely to be superficial: in fact, the equations studied by \cite{SL3Z} such as $q$-deformed heat equation, are expected to be closely related to non-stationary difference equations of the form we discuss here. The major difference, however, is that for \cite{SL3Z} many equations are very complicated. Normally a special function is defined by an equation which is much simpler than the function: for example, the Toda solution is of course an immensely complicated special function, \linebreak while Toda Hamiltonian is just a simple polynomial. Unlike this, in \cite{SL3Z} and many related works, the defining \textit{operators belong to the same class of functions as the solutions}. The same complexity is reported in \cite{Shiraishi2}. One of the goals of the present paper was to show that it is \textit{not necessarily so}: we demonstrate that for relatively non-trivial example (XXZ integrable model here) the non-stationary equation remains elementary.

Another related subject is polynomial representations of Double Affine Hecke Algebras (DAHA). \linebreak Important role in this theory is played by $SL(2,{\mathbb Z})$ automorphisms which -- without any strict necessity to introduce roots of unity -- act on polynomials in variables $x_1, x_2, \ldots$ and involve operators of the form

\begin{align}
{\widehat \gamma} \ x^n = q^{\frac{n(n+1)}{2}} \ x^n, \ \ \ n \in {\mathbb Z}
\end{align}
\smallskip\\
From this perspective, above non-stationary equations look like repeated action of such automorphisms. Again, the difference equation involving two applications of ${\widehat \gamma}$ is reminiscent of identities of the form $S^2 = \mbox{id}$ this time beyond roots of unity. It is interesting to speculate that study of non-stationary, perhaps elliptic difference equations could lead to a generalization of DAHA with at least $SL(3,{\mathbb Z})$ automorphisms.

The results of this paper have a direct connection to supersymmetric gauge theory with 8 supercharges in 5 spacetime dimensions where an interesting observable -- expectation value of a surface defect \cite{SD} -- is anticipated to satisfy a non-stationary difference equation. By construction, function $\Psi(\Lambda,x)$ is precisely the expectation value of a surface defect in 5-dimensional $SU(2)$ theory with 4 fundamental matter multiplets whose (exponentiated) masses are related to $T_1,T_2,T_3,T_4$. Eq. (\ref{NonStatEq}) is the non-stationary difference equation satisfied by this observable, and therefore 5-dimensional generalization of the results of \cite{SD-KZ}. It would be interesting to further generalize these results to 6 dimensions, obtaining some non-stationary elliptic equations and perhaps ultimately the non-stationary extension of the DELL system \cite{DELL}.

From the point of view of algebraic geometry, one feature of our construction that may be interesting is that we use representation theory of $q$-Virasoro algebra instead of quantum affine algebras, accordingly, the algebraic varieties relevant for us here are different from many varieties considered in the literature, particularly from Laumon spaces \cite{AG-Negu}. This parallels the fact in physics that there are several equivalent constructions of surface defects: current paper has to do with Higgsing defects from $SU(2) \times SU(2)$ \cite{Higgsing} while Laumon spaces are more closely related to monodromy defects \cite{DefectsAffine} and ramified instantons \cite{SD}. \linebreak It would be also interesting to understand the relation, if any, between non-stationary difference equations of the type considered here and quantum difference equations for Nakajima quiver varieties \cite{KTheoryEquation}.

Yet another direction where non-stationary difference equations may be useful is for refined topological string theory, where the same special function $\Psi(\Lambda,x)$ has the meaning of open string partition function in the presence of internal $D3$-brane (see e.g. \cite{RefTopBranes} section 5) where parameters $q$ and $t$ are both generic. \linebreak The difference equation that we obtained may be useful to study refined topological string on toric Calabi-Yau manifolds generalizing local ${\mathbb P}^1 \times {\mathbb P}^1$ with extra K\"ahler parameters corresponding to $T_1, T_2, T_3, T_4$. \linebreak One may expect generalization of our equation to other geometries, including the case of local ${\mathbb P}^2$ which is well-known for not having a direct gauge theory description. It would be interesting to explore potential benefit of non-stationary equations to the study of non-perturbative partition functions in the context of TS-ST correspondence \cite{TS-ST}.

In conformal field theory, equations satisfied by conformal blocks with degenerate insertions go back to \cite{BPZ} what is known as Belavin-Polyakov-Zamolodchikov equations. In \cite{BPZ}, these equations have been derived using existence of null-states in appropriate representations of the Virasoro algebra. Equation (\ref{NonStatEq}) provides a $q$-deformation of BPZ equations for 5-point conformal blocks with a degenerate insertion. It would be interesting to establish direct contact between equations like (\ref{NonStatEq}) and null-states of $q$-Virasoro algebra.

Last but not least, the equations considered in present paper can be interpreted as Ward identities for generalized matrix models. This is because conformal blocks in Dotsenko-Fateev representation are \linebreak naturally identified \cite{DFcomb} with $q,t$-deformed matrix models or $\beta$-ensembles \cite{qtBlock} of multi-Penner or Selberg type. \linebreak For the simplest $q,t$-deformed matrix models, Ward identities can be explicitly derived and take form of $q$-Virasoro constraints \cite{qVirConstraints} however for more general models such as the multi-Penner needed for the \linebreak description of conformal blocks, the equations get more complicated and not fully explicit \cite{qWard}. The results of current paper seem to indicate that Ward identities can be made fully explicit at the cost of becoming "exponentiated", involving infinite products in the coordinate part and ${\widehat \gamma}$ operators in the momentum part. \linebreak It would be interesting to derive the equation (\ref{NonStatEq}) as Ward identity (Stokes theorem) for the integral (\ref{qtConfBlock}).

\section{Acknowledgments}

We would like to thank M. Aganagic, S. Arthamonov, L. Chekhov, P. Etingof, A. Grekov, P. Koroteev, \linebreak C. Kozçaz, A. Marshakov, A. Mironov, A. Morozov, N. Nekrasov and A. Smirnov for valuable discussions \linebreak related to different parts of this project, and especially R. Kashaev for illuminating explanations about quantum-mechanical treatment of finite-difference problems. The computations for this work have been done using the computer algebra system MAPLE maplesoft.com/products/Maple/. The research of S.S. was supported in part by the RFBR grant 19-02-00815.


\begin{thebibliography}{11}

\bibitem{Virasoro} M.Virasoro, \emph{Subsidiary conditions and ghosts in dual-resonance models}, Phys. Rev. D. 1-10 (1970) 2933–2936, doi.org/10.1103/PhysRevD.1.2933

\bibitem{VirasoroGF} I.Gel'fand, D.Fuchs, \emph{The cohomologies of the lie algebra of the vector fields in a circle}, Funct. Anal. Appl., 2-4 (1968) 342–343, doi.org/10.1007/BF01075687

\bibitem{BPZ} A.Belavin, A.Polyakov, A.Zamolodchikov, \emph{Infinite conformal symmetry in two-dimensional quantum field theory}, Nucl. Phys. B. 241-2 (1984) 333–380,  doi.org/10.1016/0550-3213(84)90052-X

\bibitem{Polyakov} A.Polyakov, \emph{Quantum geometry of bosonic strings}, Phys. Lett. B 103-3 (1981) 207-210, doi.org/10.1016/0370-2693(81)90743-7

\bibitem{ModuliW} E.Witten, \emph{Two-dimensional gravity and intersection theory on moduli space}, Surveys in differential geometry 1 (1990) 243-310, doi.org/10.4310/SDG.1990.v1.n1.a5

\bibitem{ModuliK} M.Kontsevich, \emph{Intersection theory on the moduli space of curves and the matrix Airy function}, Comm. Math. Phys. 147-1 (1992) 1–23, doi.org/10.1007/BF02099526

\bibitem{MO} D.Maulik, A.Okounkov, \emph{Quantum Groups and Quantum Cohomology}, Asterisque 408 (2019) 212, doi.org/10.24033/ast.1074, arxiv.org/abs/1211.1287

\bibitem{DVV} R.Dijkgraaf, H.Verlinde, E.Verlinde, \emph{Loop equations and Virasoro constraints in non-perturbative two-dimensional quantum gravity}, Nucl. Phys. B 348-3(1991)435-456, doi.org/10.1016/0550-3213(91)90199-8

\bibitem{AGT} L.Alday, D.Gaiotto, Y.Tachikawa, \emph{Liouville Correlation Functions from Four-dimensional Gauge Theories}, Lett.Math.Phys. 91 (2010) 167-197, doi.org/10.1007/s11005-010-0369-5, arxiv.org/abs/0906.3219

\bibitem{VirCoproduct} G.Moore, N.Seiberg, \emph{Classical and quantum conformal field theory}, Comm. Math. Phys. 123 (1989)177–254, doi.org/10.1007/BF01238857

\bibitem{EK} P.Etingof, A.Kirillov Jr., \emph{A unified representation-theoretic approach to special functions},  Funct. Anal. Appl. 28:1 (1994) 73–76, doi.org/10.1007/BF01079011, arxiv.org/abs/hep-th/9312101

\bibitem{SF-Hyper} I.Gelfand, M.Kapranov, A.Zelevinsky, \emph{Generalized Euler integrals and A-hypergeometric functions}, Advances in Mathematics 84-2 (1990) 255-271, doi.org/10.1016/0001-8708(90)90048-R

\bibitem{SF-Gamma} S.Ruijsenaars, \emph{On Barnes' multiple zeta and gamma functions}, Advances in Mathematics 156-1 (2000) 107-132, doi.org/10.1006/aima.2000.1946

\bibitem{SF-Theta} B.Dubrovin, \emph{Theta functions and non-linear equations}, Russ. Math. Surv. 36 (1981) 11, doi.org/10.1070/RM1981v036n02ABEH002596

\bibitem{SF-Vile} N.Vilenkin, \emph{Special functions and the theory of group representations}, Amer. Math. Soc., Providence, 1968

\bibitem{SF-Wign} E.Wigner, J.Talman, \emph{Special Functions. A Group Theoretic Approach}, American Journal of Physics 37 (1969) 1073, doi.org/10.1119/1.1975207

\bibitem{EK-2} P.Etingof, A.Kirillov Jr., \emph{Macdonald’s polynomials and representations of quantum groups}, Math. Res. Let. 1 (1994) 279–296, doi.org/10.4310/MRL.1994.v1.n3.a1, arxiv.org/abs/hep-th/9312103

\bibitem{GenMac-1} P.Etingof, A.Kirillov Jr., \emph{On the affine analogue of Jack's and Macdonald's polynomials}, Duke Math. J. 78-2 (1995) 229-256, doi.org/10.1215/S0012-7094-95-07810-7, arxiv.org/abs/hep-th/9403168

\bibitem{SF-New-1} A.Alexandrov, A.Mironov, A.Morozov, \emph{Partition Functions of Matrix Models as the First Special Functions of String Theory I. Finite Size Hermitean 1-Matrix Model}, Int.J.Mod.Phys.A 19 (2004) 4127-4165, doi.org/10.1142/S0217751X04018245, arxiv.org/abs/hep-th/0310113; \\ A.Alexandrov, A.Mironov, A.Morozov, P.Putrov, \emph{Partition Functions of Matrix Models as the First Special Functions of String Theory. II. Kontsevich Model}, Int.J.Mod.Phys.A 24 (2009) 4939-4998, doi.org/10.1142/S0217751X09046278, arxiv.org/abs/0811.2825

\bibitem{SF-New-2} A.Mironov, A.Morozov, \emph{The Power of Nekrasov Functions}, Phys.Lett.B 680 (2009) 188-194, doi.org/10.1016/j.physletb.2009.08.061, arxiv.org/abs/0908.2190

\bibitem{DF} V.Dotsenko, V.Fateev, \emph{Conformal Algebra and Multipoint Correlation Functions in Two-Dimensional Statistical Models}, Nucl.Phys. B240 (1984) 312, doi.org/10.1016/0550-3213(84)90269-4

\bibitem{DFcomb} A.Mironov, A.Morozov, Sh.Shakirov, \emph{Conformal blocks as Dotsenko-Fateev integral discriminants}, Int. J. Mod. Phys. A 25-16 (2010) 3173-3207, doi.org/10.1142/S0217751X10049141, arxiv.org/abs/1001.0563

\bibitem{DFloop} A.Mironov, A.Morozov, Sh.Shakirov, \emph{On "Dotsenko-Fateev" representation of the toric conformal blocks}, J.Phys.A 44(2011) 085401, doi.org/10.1088/1751-8113/44/8/085401, arxiv.org/abs/1010.1734

\bibitem{OrdinaryW} A.Zamolodchikov, \emph{Infinite additional symmetries in two-dimensional conformal quantum field theory}, Theor. Math. Phys., 65-3 (1985) 1205–1213, doi.org/10.1007/BF01036128

\bibitem{qtVir} J.Shiraishi, H.Kubo, H.Awata, S.Odake, \emph{A Quantum Deformation of the Virasoro Algebra and the Macdonald Symmetric Functions}, Lett.Math.Phys. 38 (1996) 33, doi.org/10.1007/BF00398297, arxiv.org/abs/q-alg/9507034

\bibitem{qtVirReview} H.Awata, H.Kubo, S.Odake, J.Shiraishi, \emph{Virasoro-type Symmetries in Solvable Models}, Lett.Math.Phys. 38 (1996) 33, doi.org/10.1007/BF00398297, arxiv.org/abs/hep-th/9612233

\bibitem{qtBlock} H.Awata, Y.Yamada, \emph{Five-Dimensional AGT Relation and the Deformed $\beta$-Ensemble}, Prog. Theor. Phys. 124-2 (2010) 227–262, doi.org/10.1143/PTP.124.227, arxiv.org/abs/1004.5122

\bibitem{FK} L.Faddeev, R.Kashaev, \emph{Quantum Dilogarithm}, Mod.Phys.Lett. A 9 (1994) 427–434, doi.org/10.1142/S0217732394000447, arXiv:hep-th/9310070

\bibitem{NS} N.Nekrasov, S.Shatashvili, \emph{Quantization of Integrable Systems and Four Dimensional Gauge Theories}, XVI International Congress on Mathematical Physics 265-289, World Sci. Publ., 2010, arxiv.org/abs/0908.4052

\bibitem{AG-Negu} A.Neguţ, \emph{Affine Laumon spaces and integrable systems}, arxiv.org/abs/1112.1756

\bibitem{statDefect} M.Bullimore, H.-C.Kim, P.Koroteev, \emph{Defects and Quantum Seiberg-Witten Geometry}, JHEP 05 (2015) 095, doi.org/10.1007/JHEP05(2015)095, arxiv.org/abs/1412.6081

\bibitem{Shiraishi1} J.Shiraishi, \emph{Affine Screening Operators, Affine Laumon Spaces, and Conjectures Concerning Non-Stationary Ruijsenaars Functions}, J. Int. Sys., 4-1 (2019) xyz010, doi.org/10.1093/integr/xyz010, arxiv.org/abs/1903.07495

\bibitem{Shiraishi2} E.Langmann, M.Noumi, J.Shiraishi, \emph{Basic Properties of Non-Stationary Ruijsenaars Functions}, SIGMA 16 (2020) 105, doi.org/10.3842/SIGMA.2020.105, arxiv.org/abs/2006.07171

\bibitem{qtVirPoisson} E.Frenkel, N.Reshetikhin, \emph{Quantum Affine Algebras and Deformations of the Virasoro and ${\cal W}$-algebras}, Comm. Math. Phys. 178 (1996) 237–264, doi.org/10.1007/BF02104917, arxiv.org/abs/q-alg/9505025

\bibitem{qtVirDCA} E.Frenkel, N.Reshetikhin, \emph{Towards Deformed Chiral Algebras}, XXIth International Colloquium on Group Theoretical Methods in Physics, Goslar, 1996, arxiv.org/abs/q-alg/9706023

\bibitem{qtAGT} H.Awata, Y.Yamada, \emph{Five-dimensional AGT conjecture and the deformed Virasoro algebra}, JHEP 1001 (2010) 125, doi.org/10.1007/JHEP01(2010)125, arxiv.org/abs/0910.4431

\bibitem{5dKtheory} E.Carlsson, N.Nekrasov, A.Okounkov, \emph{Five dimensional gauge theories and vertex operators}, Mosc. Math. J.  14-1 (2014) 39-61, doi.org/10.17323/1609-4514-2014-14-1-39-61, arxiv.org/abs/1308.2465

\bibitem{Localization} M.Atiyah, R.Bott, \emph{The moment map and equivariant cohomology}, Topology 23-1 (1984) 1-28, doi.org/10.1016/0040-9383(84)90021-1

\bibitem{qAGT-2} H.Awata, B.Feigin, A.Hoshino, M.Kanai, J.Shiraishi, S.Yanagida, \emph{Notes on Ding-Iohara algebra and AGT conjecture}, arxiv.org/abs/1106.4088

\bibitem{qAGT-3} M.Taki, \emph{On AGT-${\cal W}$ Conjecture and $q$-Deformed ${\cal W}$-Algebra}, arxiv.org/abs/1403.7016

\bibitem{qAGT-Proof-1} M.Aganagic, Sh.Shakirov, \emph{Gauge/Vortex duality and AGT}, New Dualities of Supersymmetric Gauge Theories, Mathematical Physics Studies, Springer, 2016, doi.org/10.1007/978-3-319-18769-3\_13, arxiv.org/abs/1412.7132

\bibitem{qAGT-Proof-2} A.Mironov, A.Morozov, Sh.Shakirov, A.Smirnov, \emph{Proving AGT conjecture as HS duality: extension to five dimensions}, Nucl. Phys. B 855 (2011) 128-151, doi.org/10.1016/j.nuclphysb.2011.09.021, arxiv.org/abs/1105.0948

\bibitem{qAGT-Proof-3} M.Fukuda, Y.Ohkubo, J.Shiraishi, \emph{Generalized Macdonald Functions on Fock Tensor Spaces and Duality Formula for Changing Preferred Direction}, Comm. Math. Phys. 380-1 (2020) 1-70, doi.org/10.1007/s00220-020-03872-4, arxiv.org/abs/1903.05905

\bibitem{qAGT-Proof-4} A.Neguţ, \emph{The $q$-AGT–${\cal W}$ Relations Via Shuffle Algebras}, Comm. Math. Phys. 358 (2018) 101–170, doi.org/10.1007/s00220-018-3102-3, arxiv.org/abs/1608.08613

\bibitem{5dInst} A.Lawrence, N.Nekrasov, \emph{Instanton sums and five-dimensional gauge theories}, Nucl.Phys. B513 (1998) 239-265, doi.org/10.1016/S0550-3213(97)00694-9, arxiv.org/abs/hep-th/9706025

\bibitem{Convergence} G.Felder, M.M\"uller-Lennert, \emph{Analyticity of Nekrasov Partition Functions}, Comm. Math. Phys. 364 (2018) 683–718, doi.org/10.1007/s00220-018-3270-1, arxiv.org/abs/1709.05232

\bibitem{DegenFieldReview} A.Marshakov, A.Mironov, A.Morozov, \emph{On AGT Relations with Surface Operator Insertion and Stationary Limit of Beta-Ensembles}, J. Geom. Phys. 61 (2011) 1203-1222, doi.org/10.1016/j.geomphys.2011.01.012, arxiv.org/abs/1011.4491

\bibitem{Higgsing} D.Gaiotto, H.-C.Kim, \emph{Surface defects and instanton partition functions}, JHEP 10 (2016) 012, doi.org/10.1007/JHEP10(2016)012, arxiv.org/abs/1412.2781

\bibitem{SW-Integ-Mir} A.Mironov, \emph{Seiberg-Witten theory and duality in integrable systems}, Proc. XXXIV PNPI Winter School, arxiv.org/abs/hep-th/0011093

\bibitem{Macd} I.Macdonald, \emph{Symmetric functions and Hall polynomials}, 2nd ed., Oxford Mathematical Monographs, The Clarendon Press, Oxford University Press, New York, 1995

\bibitem{RS} S.Ruijsenaars, H.Schneider, \emph{A new class of integrable systems and its relation to solitons}, doi.org/10.1016/0003-4916(86)90097-7

\bibitem{Gasper-Rahman} G.Gasper, M.Rahman, \emph{Basic Hypergeometric Series}, Encyclopedia of Mathematics and its Applications,  35, Cambridge University Press, Cambridge, 1990

\bibitem{qSaalschutzJackson} F.Jackson, \emph{Transformations of $q$-series}, Messenger of Math. 39 (1910) 145-151

\bibitem{CutAndJoin1} I.Goulden, D.Jackson, \emph{Transitive factorizations into transpositions and holomorphic mappings on the sphere}, Proc. Amer. Math. Soc., 125 (1997) 51-60, jstor.org/stable/2161793

\bibitem{CutAndJoin2} M.Kazarian, \emph{KP hierarchy for Hodge integrals}, Adv. Math. 221-1 (2009) 1-21, doi.org/10.1016/j.aim.2008.10.017, arxiv.org/abs/0809.3263

\bibitem{CutAndJoin3} A.Alexandrov, \emph{Cut-and-Join operator representation for Kontsevich-Witten tau-function}, Mod. Phys. Lett. A 26 (2011) 2193-2199, doi.org/10.1142/S0217732311036607, arxiv.org/abs/1009.4887

\bibitem{Wrep1} A.Morozov, Sh.Shakirov, \emph{Generation of matrix models by ${\cal W}$-operators}, JHEP 0904 (2009) 064, doi.org/10.1007/JHEP04(2009)064, arxiv.org/abs/0902.2627

\bibitem{Wrep2} A.Morozov, Sh.Shakirov, \emph{On ${\cal W}$-representations of $\beta$- and $q,t$-deformed matrix models}, Phys.Lett. B792 (2019) 205-213, doi.org/10.1016/j.physletb.2019.03.047, arxiv.org/abs/1901.02811

\bibitem{Wrep3} A.Mironov, V.Mishnyakov, A.Morozov, \emph{Non-Abelian ${\cal W}$-representation for GKM}, Phys.Lett. B 823 (2021) 136-721, doi.org/10.1016/j.physletb.2021.136721, arxiv.org/abs/2107.02210

\bibitem{KashaevOperators} R.Kashaev, \emph{Discrete Liouville equation and Teichmüller theory}, Handbook of Teichm\"uller theory. III, IRMA Lect. Math. Theor. Phys. 17 (2012) 821–851, doi.org/10.4171/103-1/16, arxiv.org/abs/0810.4352

\bibitem{inner} A.Kirillov Jr., \emph{On inner product in modular tensor categories}, J. Amer. Math. Soc. 9 (1996) 1135-1169, doi.org/10.1090/S0894-0347-96-00210-X, arxiv.org/abs/q-alg/9508017

\bibitem{SL3Z} G.Felder, A.Varchenko, \emph{$q$-deformed KZB heat equation: completeness, modular properties and $SL(3,{\mathbb Z})$}, Adv. Math. 171(2002)228-275, doi.or/10.1006/aima.2002.2080, arxiv.org/abs/math/0110081

\bibitem{SD} S.Gukov, E.Witten, \emph{Gauge Theory, Ramification, And The Geometric Langlands Program}, Current developments in mathematics, Int. Press, Somerville, 2008, 35–180, arxiv.org/abs/hep-th/0612073

\bibitem{SD-KZ} N.Nekrasov, A.Tsymbaliuk, \emph{Surface defects in gauge theory and KZ equation}, arxiv.org/abs/2103.12611

\bibitem{DELL} P.Koroteev, Sh.Shakirov, \emph{The Quantum DELL System}, Lett. Math. Phys. 110 (2020) 969–999, doi.org/10.1007/s11005-019-01247-y, arxiv.org/abs/1906.10354

\bibitem{DefectsAffine} L.Alday, Y.Tachikawa, \emph{Affine $SL(2)$ Conformal Blocks from 4d Gauge Theories}, Lett. Math. Phys., 94 (2010) 87–114, doi.org/10.1007/s11005-010-0422-4, arxiv.org/abs/1005.4469

\bibitem{KTheoryEquation} A.Okounkov, A.Smirnov, \emph{Quantum difference equation for Nakajima varieties}, \\ arxiv.org/abs/1602.09007

\bibitem{RefTopBranes} C.Kozçaz, Sh.Shakirov, C.Vafa, W.Yan, \emph{Refined Topological Branes}, Comm. Math. Phys. 385 (2021) 937–961, doi.org/10.1007/s00220-020-03883-1,  arxiv.org/abs/1805.00993

\bibitem{TS-ST} A.Grassi, Y.Hatsuda, M.Marino, \emph{Topological Strings from Quantum Mechanics}, Ann. Henri Poincaré 17 (2016) 3177–3235, doi.org/10.1007/s00023-016-0479-4, arxiv.org/abs/1410.3382

\bibitem{qVirConstraints} R.Lodin, A.Popolitov, Sh.Shakirov, M.Zabzine, \emph{Solving $q$-Virasoro constraints}, Lett. Math. Phys. 110 (2020) 179–210, doi.org/10.1007/s11005-019-01216-5, arxiv.org/abs/1810.00761

\bibitem{qWard} Y.Zenkevich, \emph{Quantum spectral curve for $(q, t)$-matrix model}, Lett. Math. Phys. 108 (2018) 413–424, doi.org/10.1007/s11005-017-1015-2, arxiv.org/abs/1507.00519

\end{thebibliography}
\end{document}